# TRACY–WIDOM LIMIT FOR THE LARGEST EIGENVALUE OF A LARGE CLASS OF COMPLEX SAMPLE COVARIANCE MATRICES[1]

By Noureddine El Karoui

*University of California, Berkeley*

We consider the asymptotic fluctuation behavior of the largest eigenvalue of certain sample covariance matrices in the asymptotic regime where both dimensions of the corresponding data matrix go to infinity. More precisely, let $X$ be an $n \times p$ matrix, and let its rows be i.i.d. complex normal vectors with mean 0 and covariance $\Sigma_p$. We show that for a large class of covariance matrices $\Sigma_p$, the largest eigenvalue of $X^*X$ is asymptotically distributed (after recentering and rescaling) as the Tracy–Widom distribution that appears in the study of the Gaussian unitary ensemble. We give explicit formulas for the centering and scaling sequences that are easy to implement and involve only the spectral distribution of the population covariance, $n$ and $p$.

The main theorem applies to a number of covariance models found in applications. For example, well-behaved Toeplitz matrices as well as covariance matrices whose spectral distribution is a sum of atoms (under some conditions on the mass of the atoms) are among the models the theorem can handle. Generalizations of the theorem to certain spiked versions of our models and a.s. results about the largest eigenvalue are given. We also discuss a simple corollary that does not require normality of the entries of the data matrix and some consequences for applications in multivariate statistics.

**1. Introduction.** Sample covariance matrices are a fundamental tool of multivariate statistics. In the classical setting, one starts with an $n \times p$ data matrix $X$ and studies asymptotic properties of $S = (X - \bar{X})'(X - \bar{X})/(n-1)$ when $p$ is fixed and $n$ grows to infinity. The classic paper [1] answered most

Received May 2005; revised May 2006.

[1]Supported in part by NSF Grant DMS-06-05169.

of the relevant questions concerning the eigenvalues of $S$ in the setting where the rows of $X$ are i.i.d. $\mathcal{N}(M,\Sigma)$. It was shown in [1] that, from an eigenvalue point of view, $S$ was a good estimator of $\Sigma$. A thorough account of the classical case can be found in [2], Chapters 11 and 13.

Nowadays, statisticians are working with datasets of increasingly larger size and the practical relevance of the assumption that $p$ is fixed and $n$ goes to infinity is often doubtful. It might also be counterproductive in applications. A significant effort has therefore been made recently to try to understand the asymptotic behavior of certain classical tools in multivariate analysis, such as the largest eigenvalue of $S$, in the setting where $p$ and $n$ both grow to infinity.

Large-dimensional sample covariance matrices are also of interest in other fields than statistics. Matrices $X$ whose entries are complex-valued and their singular values are also of interest in different fields of applications and in particular in communications engineering. They are objects of great interest, for instance, in wireless communications (see, e.g., [32] and [38]).

In the rest of the paper, the data matrix will always be called $X$. The eigenvalues of $X^*X$ will be denoted $l_1 \geq l_2 \geq \cdots \geq l_p$. The population eigenvalues, that is, the eigenvalues of $\Sigma_p$ will be called $\lambda_1 \geq \lambda_2 \geq \cdots \geq \lambda_p$.

To situate our paper in the current literature, let us recall a few results that have been recently obtained. When the true covariance matrix is Id and the entries of $X$ are either standard complex or standard real normal distributed, results in [14, 24, 25] and [12] showed that

$$\text{if } n \text{ and } p \to \infty, \qquad \frac{l_1(X^*X)-\mu_{n,p}}{\sigma_{n,p}} \Rightarrow \mathrm{TW},$$

where $\mu_{n,p}$ and $\sigma_{n,p}$ are explicit sequences (which do not depend on whether the real or the complex case is under consideration), and the limiting law is a Tracy–Widom distribution. When the entries are standard complex normal, the limiting law is the Tracy–Widom distribution appearing in the study of the Gaussian unitary ensemble (see [35]). When they are real normal, it is the one corresponding to the Gaussian orthogonal ensemble (see [36]).

More recently, the paper [6] looked at finite-dimensional perturbations of the Id covariance matrix. The authors considered so-called "spiked" covariance models, advocated in [25], where a finite number—$k$—of eigenvalues are changed to a value different from 1 and the remaining $p-k$ eigenvalues are fixed at 1. They discovered a very interesting phase transition phenomenon, with the behavior of $l_1$ changing drastically depending on how far away $\lambda_1$, the largest eigenvalue of $\Sigma_p$, is from the bulk of the spectrum of $\Sigma_p$. In their case, this bulk was of course concentrated at 1.

In the course of their analysis, they develop powerful tools to analyze the problem. In particular, their Proposition 2.1 (for which they also give credit to K. Johansson), and the subsequent remarks are finite dimensional and valid whatever the true covariance structure. We exploit in this paper the powerful representations obtained in [6] to handle a much more general class of covariance matrices than finite perturbations of the Id matrix.

The motivations for doing so are many. From a theoretical standpoint, it is somewhat unclear at this point what features, if any, of the covariance structure of the random variables are responsible for the appearance of Tracy–Widom laws. One might ask, for instance, if it is the fact that the bulk of the true eigenvalues is exactly concentrated at one point. We will show that intuitively what seems needed is a weaker condition, the clumping of a fraction of eigenvalues close to the largest one.

From an applications standpoint, many covariances appearing in different fields of science are not finite-dimensional perturbations of the Id matrix. Block-diagonal covariance matrices are of particular interest since they are accepted models for, say, the correlation of genes in microarray analysis (a topic of intense statistical research at the time being), or the correlation of the returns of stocks of companies in financial applications. Covariances that are sums of atoms, for example, $a\%$ of the variables have variance $\lambda_1$ and $1-a\%$ have variance $\lambda_2$, are also of interest, especially in light of Theorem 1.1b in [6]. We will come back to this in Section 4. In other respects, covariance matrices that are also Toeplitz matrices are very natural in the analysis of time-series data, since the covariance structure of a stationary time-series is a Toeplitz matrix.

Before we state our main theorem, we need to introduce some terminology and set some notation. We will be working with $n \times p$ matrices $X$, whose rows $\{X_k\}_{k=1,\dots,n}$ are i.i.d. $\mathcal{N}_{\mathbb{C}}(0,\Sigma_p)$. By definition, this means that $X_k = Y_k + iZ_k$, where $Y_k$ and $Z_k$ are independent (real) $\mathcal{N}(0,\Sigma_p/2)$. The matrix $W = X^*X$ is then called a complex Wishart matrix, with $n$ degrees of freedom and covariance $\Sigma_p$. It will be abbreviated $W_{\mathbb{C}}(\Sigma_p, n)$.

We will call the eigenvalues of $\Sigma_p$ $\lambda_i$, with $\lambda_1 \geq \lambda_2 \geq \cdots \geq \lambda_p$. The eigenvalues of $W_{\mathbb{C}}(\Sigma_p, n)$ will be denoted $l_i$, with the same ordering convention, that is, $l_1$ is the largest eigenvalue of $X^*X$.

It is well known in statistics that if $X_k$ are i.i.d. $\mathcal{N}_{\mathbb{C}}(M,\Sigma_p)$, then $(X-\bar{X})^*(X-\bar{X})$ is $W_{\mathbb{C}}(\Sigma_p, n-1)$. For this reason, we will always assume that the $X_k$'s are $\mathcal{N}_{\mathbb{C}}(0,\Sigma_p)$.

We are now ready to state the main theorem.

THEOREM 1. *Let us consider complex Wishart matrices $W_{\mathbb{C}}(\Sigma_p, n)$. Let $\lambda_1$ be the largest eigenvalue of $\Sigma_p$ and let $\lambda_p$ be the smallest one. Let $H_p$ be the spectral distribution of $\Sigma_p$. Let $c$ be the unique solution in $[0, 1/\lambda_1(\Sigma_p))$ of the equation*

$$c = c(\Sigma_p, n, p), c \in [0, 1/\lambda_1(\Sigma_p)) : \int \left(\frac{\lambda c}{1-\lambda c}\right)^2 dH_p(\lambda) = \frac{n}{p}. \tag{1}$$

*We assume that* $n/p \geq 1$ *is uniformly bounded,* $\limsup \lambda_1 < \infty$, $\liminf \lambda_p > 0$ *and* $\limsup \lambda_1 c < 1$. *We denote by* $\mathcal{G}$ *the class of models* $\{(\Sigma_p, n, p)\}$ *for which these conditions hold. We call*

$$\mu = \frac{1}{c}\left(1 + \frac{p}{n}\int \frac{\lambda c}{1-\lambda c}\, dH_p(\lambda)\right), \tag{2}$$

$$\sigma^3 = \frac{1}{c^3}\left(1 + \frac{p}{n}\int \left(\frac{\lambda c}{1-\lambda c}\right)^3 dH_p(\lambda)\right). \tag{3}$$

*Let* $l_1$ *be the largest eigenvalue of* $W_{\mathbb{C}}(\Sigma_p, n)$, *that is,* $l_1 = l_1(X^*X)$, *where* $X$ *is an* $n \times p$ *matrix whose rows are i.i.d.* $\mathcal{N}_{\mathbb{C}}(0, \Sigma_p)$. *Then we have, as* $n$ *goes to* $\infty$,

$$\frac{l_1 - n\mu}{\sigma n^{1/3}} \quad \Longrightarrow \quad \mathrm{TW}_2.$$

*Moreover, if we denote by* $F_0$ *the cumulative distribution function of* $\mathrm{TW}_2$, *we can find* $\varepsilon > 0$ *and a continuous, nonincreasing function* $C$ *(that may depend on the models under consideration and* $\varepsilon$*) such that*

$$\forall s_0\ \exists N_0 : s \geq s_0 \text{ and } n \geq N_0 \text{ implies}$$
$$\left|P\left(\frac{l_1 - n\mu}{n^{1/3}\sigma} \leq s\right) - F_0(s)\right| \leq \frac{C(s_0)e^{-\varepsilon s/2}}{n^{1/3}}.$$

Using these results, their proofs, and a little bit more work, we can prove the following corollaries:

COROLLARY 1. *In the setting of Theorem 1, if* $\{(\Sigma_p, n, p)\}$ *is in* $\mathcal{G}$, *we have*

$$\frac{l_1}{n} - \mu \to 0 \qquad a.s.$$

COROLLARY 2. *In the setting of Theorem 1, if* $\{(\Sigma_p, n, p)\}$ *is in* $\mathcal{G}$, *the* $k$*-largest eigenvalues of* $X^*X$, *properly recentered and rescaled, converge to their Tracy–Widom counterpart.*

Before we proceed, let us remind the reader that the cumulative distribution function of $\mathrm{TW}_2$ is known. After introducing the intermediary function $q$ defined by

$$q''(x) = xq(x) + 2q^3(x),$$
$$q(x) \sim \mathrm{Ai}(x) \qquad \text{as } x \to \infty,$$

$F_0$ satisfies (see [35])

$$F_0(s) = \exp\left(-\int_s^\infty (x-s)q^2(x)\, dx\right).$$

We will discuss in greater detail the potential usage of the theorem in Section 4, but we want to highlight sufficient conditions under which it applies and give a few examples before we give the proof. More examples, additional results concerning spiked versions of matrices in $\mathcal{G}$ and a remark about the fact that the bias of $l_1$ is (in some cases) independent of the distributional assumptions made on the entries will be found in Section 4.

COROLLARY 3 (Sufficient conditions). *When the following five conditions are all satisfied, the theorem applies:*

1. *$n/p$ remains bounded and $n \geq p$;*
2. *$H_p \Rightarrow H_\infty$, in the usual weak convergence sense;*
3. *$\lambda_1(\Sigma_p) \to \lambda_1(\infty) = \sup \operatorname{support} H_\infty < \infty$;*
4. *$\lambda_p(\Sigma_p) \to \lambda_\infty(\infty) = \inf \operatorname{support} H_\infty > 0$;*
5. *$H_\infty$ has a density $h_\infty(\lambda)$ in a (left) neighborhood of $\lambda_1(\infty)$, and in this neighborhood, $h_\infty(\lambda) \geq B(\lambda_1(\infty) - \lambda)$ for some $B > 0$.*

As a consequence we see that the result applies to:

- Symmetric Toeplitz matrices—with parameters $a_0, a_1, \ldots$—for which $\sum k|a_k| < \infty$, the function

$$a : a(\omega) = a_0 + 2 \sum_{k=1}^{\infty} a_k \cos(k\omega)$$

  has a derivative that changes sign only a finite number of times on $[0, 2\pi]$, and for which the distribution $F$ of $a$ does not have atoms ($F(x) = \frac{1}{2\pi} \operatorname{Leb}\{\omega \in [0, 2\pi] : f(\omega) \leq x\}$).
- Covariances that have uniformly spaced eigenvalues on an interval $[\zeta, \xi]$, as long as $\zeta > 0$ and $\xi < \infty$.

Also, as shown in Appendix A.3.1, if $H_p$ has an atom of mass $\nu(p)$ at $\lambda_1(\Sigma_p)$ and $\liminf \nu(p) > 0$, assuming that $\limsup \lambda_1 < \infty$, $n/p$ remains bounded and $\liminf \lambda_p(\Sigma_p) > 0$, the theorem holds.

Hence the Id case, which was investigated in [14, 24] and [25], is a special case of our main theorem. Also, since spiked models with a "small" spike are in $\mathcal{G}$ (see Section 4), the results of [6] showing convergence to $\mathrm{TW}_2$ are also a subcase of our main result.

**2. Framework.** As is—almost—classical for this problem, one tries to represent the marginal distribution, $P(l_1/n \leq x)$, as the determinant of $I - K_{n,p}$, where $K_{n,p}$ is a trace class operator acting on $L^2([x, \infty))$. It greatly simplifies the analysis if one is able to represent $K_{n,p}$ as the product of Hilbert–Schmidt operators, say $H_{n,p} J_{n,p}$. The problem is even more tractable if the kernels of those operators have the property that $H_{n,p}(x,y) = H_{n,p}(x+y)$, and similarly for $J_{n,p}$.

Let us mention before we proceed that we will be denoting the trace class norm of an operator $K$ by $\|K\|_1$. Its Hilbert–Schmidt norm will be denoted by $\|K\|_2$. An introduction to these concepts can be found in [30], Section VI.6 or [16], Chapter 4.

2.1. *Finite-dimensional representation of operators.* Proposition 2.1 in [6] and their remarks in (82)–(85) remarkably managed to obtain all the characteristics of the representations we wished for in the case of completely general $\Sigma_p$. Since the authors of [6] credit Johansson for the very elegant proof they present, we will call this theorem Baik–Ben Arous–Johansson–Péché. Here is what it states.

THEOREM 2 (Baik–Ben Arous–Johansson–Péché). *Let us consider an $n \times p$ matrix $X$ with rows i.i.d. $\mathcal{N}_{\mathbb{C}}(0,\Sigma_p)$. Let us assume without loss of generality that $\lambda_p(\Sigma_p) > 0$. Let $\pi_i = 1/\lambda_i$. Let $q \in \mathbb{R}$ be such that $0 < q < \pi_1$. Let $K_{n,p}$ be the operator on $L^2([s,\infty))$ with kernel*

$$K_{n,p}(x,y) = \frac{n}{(2\pi i)^2}\int_\Gamma dz \int_\Xi dw \times \exp(-xn(z-q) + yn(w-q))\frac{1}{w-z}\left(\frac{z}{w}\right)^n \prod_{k=1}^p \frac{\pi_k - w}{\pi_k - z}. \tag{4}$$

*Here $\Xi$ (resp. $\Gamma$) is a simple closed contour oriented counterclockwise and encircling 0 (resp. $\pi_1,\dots,\pi_p$). Then, if we denote by $l_1$ the largest eigenvalue of $X^*X$, we have*

$$P\left(\frac{l_1}{n} \le s\right) = \det(I - K_{n,p}|_{L^2([s,\infty))}). \tag{5}$$

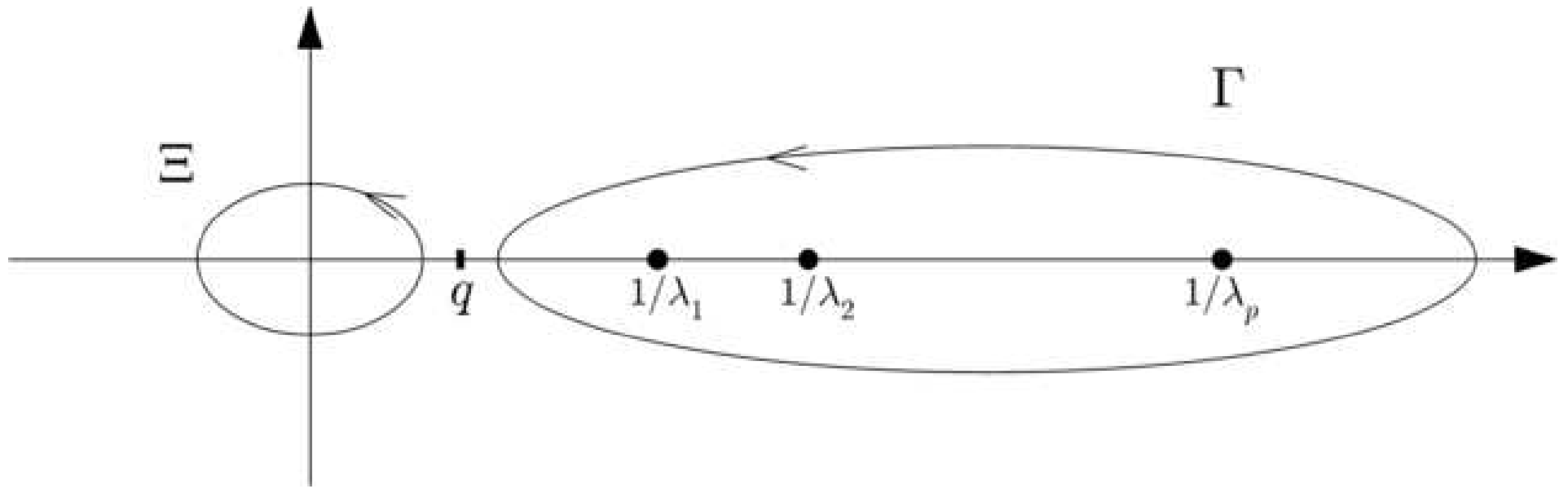

FIG. 1. *Graphical depiction of $\Gamma$ and $\Xi$.*

*Moreover, $K_{n,p}$ can be rewritten as*

$$K_{n,p}(x,y) = \int_0^\infty H_{n,p}(x+u) J_{n,p}(u+y)\, du, \tag{6}$$

*with*

$$H_{n,p}(x) = \frac{n}{2\pi} \int_\Gamma e^{-xn(z-q)} z^n \prod_{k=1}^p \frac{1}{\pi_k - z}\, dz, \tag{7}$$

$$J_{n,p}(y) = \frac{n}{2\pi} \int_\Xi e^{yn(w-q)} w^{-n} \prod_{k=1}^p (\pi_k - w)\, dw. \tag{8}$$

*Note that $\Xi$ should be strictly to the left of $\Gamma$.*

We reproduce their Figure 2 as our Figure 1 for the convenience of the reader to give a graphical representation of $\Gamma$ and $\Xi$. We refer the reader to Remark 2.1 in [6] for a discussion of the meaning of $q$. For our purposes, it will be enough to know that $q$ is essentially a free parameter that regularizes the operators we deal with.

2.2. *Recentering, rescaling and classical operator theory arguments leading to weak convergence.* Once the very important representations mentioned in (6)–(8) are obtained, the path to showing weak convergence is classical in this type of problem. One needs to find centering and scaling sequences such that the recentered and rescaled version of $K_{n,p}$ converges in trace class norm to its limit. Trace class norm plays an important role because the determinant $\det(I - \cdot)$ is continuous with respect to that norm.

Since (6) shows that $K_{n,p} = H_{n,p} J_{n,p}$, the problem reduces to showing convergence in Hilbert–Schmidt norm of $H_{n,p}$ and $J_{n,p}$ (once again properly recentered and rescaled) to their limit. This comes essentially from the fact that if $A$ and $B$ are Hilbert–Schmidt operators, $AB$ is trace class with $\|AB\|_1 \le \|A\|_2 \|B\|_2$ and some elementary algebra.

The authors of [6], in their Section 2.2, prepare the rest of their paper by doing recentering and scaling of the operators already specializing to the case of interest to them, namely finite-dimensional perturbations of the Id matrix. We do it here for general $\Sigma_p$.

Let us be more explicit now that we have explained the basic ideas. Because (5) is exact in finite dimension, one has

$$P\left(\frac{l_1}{n} \le \mu_{n,p} + \sigma_{n,p} s\right) = \det(I - K_{n,p}|_{L^2(\mu_{n,p}+\sigma_{n,p}s)}) = \det(I - S_{n,p}|_{L^2(0,\infty)})$$

and $S_{n,p}$ has kernel

$$S_{n,p}(x,y) = \sigma_{n,p} K_{n,p}(\mu_{n,p} + \sigma_{n,p}(x+s), \mu_{n,p} + \sigma_{n,p}(y+s)).$$

This is what we called earlier the recentered and rescaled operator. Because of the representation given in (6), we see that

$$S_{n,p}(x,y) = \int_0^\infty \widetilde{H}_{n,p}(x+s+u)\widetilde{J}_{n,p}(y+s+u)\,du,$$

with

$$\widetilde{H}_{n,p}(x) = \frac{n\sigma_{n,p}}{2\pi}\int_\Gamma e^{-n\sigma_{n,p}x(z-q)}e^{-n\mu_{n,p}(z-q)}z^n\prod_{k=1}^p \frac{1}{\pi_k - z}\,dz,$$

$$\widetilde{J}_{n,p}(x) = \frac{n\sigma_{n,p}}{2\pi}\int_\Xi e^{n\sigma_{n,p}x(z-q)}e^{n\mu_{n,p}(z-q)}\frac{1}{z^n}\prod_{k=1}^p(\pi_k - z)\,dz.$$

From an operator-theoretic standpoint, the three formulae above mean that

$$S_{n,p} = \widetilde{H}_{n,p}\widetilde{J}_{n,p},$$

and we can now view them as operators acting on $L^2([s,\infty))$ with kernel, that is, $\widetilde{H}_{n,p}(x,y) = \widetilde{H}_{n,p}(x+y-s)$. Now, since $\pi_k = 1/\lambda_k$, it is clear that

$$\widetilde{H}_{n,p}(x) = \frac{n\sigma_{n,p}}{2\pi}\det(\Sigma_p)\int_\Gamma e^{-n\sigma_{n,p}x(z-q)}e^{-n\mu_{n,p}(z-q)}z^n\prod_{k=1}^p\frac{1}{1-z\lambda_k}\,dz$$

and

$$\widetilde{J}_{n,p}(x) = \frac{n\sigma_{n,p}}{2\pi\det(\Sigma_p)}\int_\Xi e^{n\sigma_{n,p}x(z-q)}e^{n\mu_{n,p}(z-q)}\frac{1}{z^n}\prod_{k=1}^p(1-\lambda_k z)\,dz.$$

Being primarily interested in the product $\widetilde{H}_{n,p}\widetilde{J}_{n,p}$ and not the individual operators, we see (with [6]) that we have a little bit of choice in the operators we wish to work with. In particular, we can choose to work with $\kappa_{n,p}\widetilde{H}_{n,p}$ and $\widetilde{J}_{n,p}/\kappa_{n,p}$ for any nonzero sequence $\kappa_{n,p}$. So we can get rid of the $\det(\Sigma_p)$ term appearing in the previous display and work with

$$A_{n,p}(x) = -\frac{n\sigma_{n,p}}{2\pi i}\int_\Gamma e^{-n\sigma_{n,p}x(z-q)}e^{-n\mu_{n,p}(z-q)}z^n\prod_{k=1}^p\frac{1}{1-z\lambda_k}\,dz, \tag{9}$$

$$B_{n,p}(x) = \frac{n\sigma_{n,p}}{2\pi i}\int_\Xi e^{n\sigma_{n,p}x(z-q)}e^{n\mu_{n,p}(z-q)}\frac{1}{z^n}\prod_{k=1}^p(1-\lambda_k z)\,dz. \tag{10}$$

We now have $S_{n,p} = A_{n,p}B_{n,p}$ and $A_{n,p}$ and $B_{n,p}$ are operators on $L^2([s,\infty))$ with kernels $A_{n,p}(x,y) = A_{n,p}(x+y-s)$ and similarly for $B_{n,p}$. Since we are aiming to show convergence to $\mathrm{TW}_2$, the Airy function will play a central role in our analysis. We will denote it by Ai. Showing weak convergence of $l_1(X^*X)$ to the Tracy–Widom law reduces to finding $\kappa_{n,p}$ and "good" $A_\infty$ and $B_\infty$ such that $\|\kappa_{n,p}A_{n,p} - A_\infty\|_2 \to 0$ and $\|B_{n,p}/\kappa_{n,p} - B_\infty\|_2 \to 0$. Since, if we view the operators as acting on $L^2([s,\infty))$, $A_{n,p}(x,y) = A_{n,p}(x+y-s)$ and similarly for the Airy operator, $\mathrm{Ai}(x,y) = \mathrm{Ai}(x+y-s)$, this will essentially amount to just showing that $\kappa_{n,p}A_{n,p}(x) - A_\infty(x) \to 0$ pointwise, $A_\infty$ being a simple modification of the Airy function, and that both functions go to 0 fast enough (e.g., faster than $e^{-bx}$ for some $b > 0$) at $\infty$.

The operator-theoretic arguments used to prove Theorem 2 have considerably simplified the problem, at least conceptually: we have moved from the problem of studying an integral in $\mathbb{R}^p$ to that of analyzing a function of one real variable. Note that this was also the case with previous studies (see [25]), where arguments from [37] and [41] (where some of the ideas behind Theorem 2 can be found) were used to reduce the complexity of the problem to the same degree.

What is left to do now is very clear. We just need to find $\mu_{n,p}$, $\sigma_{n,p}$, $\Gamma$, $\Xi$ and $\kappa_{n,p}$ such that $\kappa_{n,p}A_{n,p} - A_\infty$ and $B_{n,p}/\kappa_{n,p} - B_\infty$ go to 0 (in Hilbert–Schmidt norm) when $n, p$ go to $\infty$, for appropriate $A_\infty$ and $B_\infty$. More details on these functions will be found in Propositions 1 and 2. The next section will be devoted to doing all of this.

**3. Proof of the main result.** A point of terminology before we proceed: we will interchangeably call $A_{n,p}$ either the operator whose kernel is $A_{n,p}(x+y)$ or the corresponding function. This simplifies the notation and the exposition. If there is some ambiguity, we will say precisely if we refer to the operator, its kernel or the function that defines the kernel.

The strategy of the proof is the same as that of [6]. Loosely speaking, the functions $A_{n,p}$ and $B_{n,p}$ can be viewed as integrals depending on parameters going to $\infty$. The functions to integrate contain elements of the type $e^{nf_{n,p}(z)}$. This is a situation where one can try to use steepest descent analysis.

3.1. *Focus of the analysis.* The expression defining $A_{n,p}$ in (9) is somewhat involved, but we will concentrate mostly on

$$e^{-n\mu_{n,p}(z-q)}z^n \prod_{k=1}^{p} \frac{1}{1 - z\lambda_k},$$

which can be rewritten as

$$f(z) \triangleq \exp\left(-n\mu_{n,p}(z-q) + n\log(z) - \sum_{k=1}^{p} \log(1 - z\lambda_k)\right)$$

wherever this expression makes sense. [We use the principal branch of the log, $\log(z) = \log(|z|) + i\arg(z), -\pi < \arg(z) < \pi$.] The sum appearing in the

definition of $f$ can be rewritten as an integral against the spectral distribution of $\Sigma_p$, a distribution we call $H_p$, and we finally get

$$f(z) = -\mu_{n,p}(z - q) + \log(z) - \frac{p}{n} \int \log(1 - z\lambda)\, dH_p(\lambda). \tag{11}$$

It is clear that $f$ depends on $\Sigma_p$, $n$ and $p$ but we choose to not highlight this dependence here to avoid cumbersome notation.

3.2. *Heuristic connection with work on a.s. convergence.* Many results have been obtained concerning the almost sure (a.s.) convergence of different spectral characteristics of random covariance matrices, starting with the Marčenko–Pastur equation (see [27] and [40]). The article [3] contains a thorough review and a nice introduction to these problems.

Of particular interest to us are results concerning the behavior of the largest eigenvalue in the case of non-Id covariance. Classical ([27], equation (1.15)) and more recent results (see, e.g., [4, 31]) emphasize the role of the function

$$g_\infty(m) = -\frac{1}{m} + \frac{p}{n} \int \frac{\lambda}{1 + \lambda m}\, dH_\infty(\lambda),$$

where $H_\infty$ is the limiting spectral distribution of $\Sigma_p$, in obtaining almost sure convergence properties of $l_1(X^*X/n)$ and $l_p(X^*X/n)$ and determining the limiting spectral distribution of $X^*X/n$. In particular, the points $m$ where $g'_\infty(m) = 0$ intuitively play a crucial role in determining its support. Note that doing asymptotic analysis at fixed spectral distribution and $p/n$ would lead to considering the equivalent of $g_\infty$ where $H_p$ replaces $H_\infty$.

Now proceeding formally, we see that $f'(z) = -\mu_{n,p} + g_p(-z)$, where $g_p(m) = -\frac{1}{m} + \frac{p}{n} \int \frac{\lambda}{1+\lambda m}\, dH_p(\lambda)$, and hence

$$f''(z) = -g'_p(-z).$$

Since we are essentially interested in the points where $g'_p(z) = 0$, the heuristic tells us that for a large class of $\Sigma_p$, the critical point of interest to us is going to be a triple point of the function $f$ (i.e., a saddle point of order 2).

3.3. *Consequences*: *choice of* $c, \mu_{n,p}$ *and* $\sigma_{n,p}$. So it is now clear that the solutions of the equation

$$f''(z) = -\frac{1}{z^2} + \frac{p}{n} \int \left(\frac{\lambda}{1 - \lambda z}\right)^2 dH_p(\lambda) = 0$$

are likely candidates to play a central role in the problem.

Since we are focusing on largest eigenvalue problems, it is natural to consider for $c$ the unique solution in $[0, 1/\lambda_1(\Sigma_p))$ of this equation. In other words,

$$c = c(\Sigma_p, n, p), c \in [0, 1/\lambda_1(\Sigma_p)) : \int \left(\frac{\lambda c}{1 - \lambda c}\right)^2 dH_p(\lambda) = \frac{n}{p} \tag{12}$$

will play a crucial role in our analysis.

Note that, if $a > 0$, the function $x \mapsto ax/(1-ax)$ is continuous and (strictly) increasing on $(0, 1/a)$. Hence $h(x) = \int (\lambda x)^2/(1-\lambda x)^2\, dH_p(\lambda)$ is increasing on $(0, 1/\lambda_1(\Sigma_p))$. It is also strictly convex, as a convex combination of strictly convex functions. Since $h$ goes from 0 to $\infty$ on $[0, 1/\lambda_1(\Sigma_p))$, the equation $h(x) = r$ has exactly one solution on $[0, 1/\lambda_1(\Sigma_p))$ for all $r \in \mathbb{R}_+$. Existence and uniqueness of $c(\Sigma_p, n, p)$ are therefore proved.

For steepest descent reasons, we also naturally "require" that $f'(c) = 0$ and hence

$$\mu_{n,p} = \frac{1}{c} + \frac{p}{n}\int \frac{\lambda}{1-\lambda c}\, dH_p(\lambda) = \frac{1}{c}\left(1 + \frac{p}{n}\int \frac{\lambda c}{1-\lambda c}\, dH_p(\lambda)\right).$$

Hindsight from the analysis (see Appendix A.2) makes clear that if the arguments are to go through, we will have

$$f^{(3)}(c) = 2\sigma_{n,p}^3 n^2 = \frac{2}{c^3}\left(1 + \frac{p}{n}\int \left(\frac{\lambda c}{1-\lambda c}\right)^3 dH_p(\lambda)\right).$$

While this discussion does not show anything, it provides heuristic reasons for the not necessarily intuitive choice of the parameters $c$, $\mu$ and $\sigma$. What is left to do is to find paths $\Gamma$ and $\Xi$ on which we understand the behavior of $f(z)$ and will allow us to show convergence of $A_{n,p}$ and $B_{n,p}$ to our target functions. Note that the paths $\Gamma$ and $\Xi$ we will choose are functions of $\Sigma_p$, $n$ and $p$. The fact that $f$ is real for real $z$ as well as geometric properties of saddle points of order 2 (see [28], page 137) makes natural the choice of lines crossing the real axis at angle $\pi/3$ and $2\pi/3$ as starting points for $\Gamma$ and $\Xi$ at $c$.

3.4. *About* $\Gamma$. Because of a slight technical problem appearing in the operator convergence analysis, we will not exhibit $\Gamma$ immediately but rather a $\widetilde{\Gamma}$ which is much more natural from the point of view of the analysis of the behavior of $f$. Specifically, we will exhibit a curve $\widetilde{\Gamma}_+$ on which $f(z)$ is well understood. Then $\widetilde{\Gamma}$ will just be $\widetilde{\Gamma} = \widetilde{\Gamma}_+ \cup \overline{\widetilde{\Gamma}_+}$, where the $\overline{\phantom{x}}$ denotes complex conjugation. The problem is very graphical, so we will first show a drawing of $\widetilde{\Gamma}_+$.

We will show the following lemma:

LEMMA 1. *Under the assumptions of Theorem 1,* $\Re(f(z))$ *is decreasing for* $z \in \Gamma_1 \cup \Gamma_2 \cup \Gamma_3$ *as* $\Re(z)$ *increases. Also the length of* $\widetilde{\Gamma}_+$ *is uniformly bounded. Finally, there exists* $R_1 > 0$ *such that* $\max_{z\in\Gamma_4} \Re(f(z)) \le \Re(f(d))$, *where* $d = d(\Sigma_p, n, p) = c(1 + 2(-1 + 1/(\lambda_1 c))e^{i\pi/3})$.

$R_1$ is uniform with respect to our models: given a family of models $\{(\Sigma_p, n, p)\}_{n=1}^\infty$ in $\mathcal{G}$, we get $\bar{\alpha}_1 = \limsup \lambda_1 c$, $\bar{\gamma}^2 = \limsup n/p$, $\underline{\alpha}_\infty = \liminf \lambda_p$. $R_1$ is just a function of these parameters and not of the individual triplet $(\Sigma_p, n, p)$ we will be dealing with.

Precise definitions of $\Gamma_i$'s will be given as they arise in the analysis. In particular, that of $\Gamma_2$ requires a significant amount of notation and we choose to postpone it in the interest of clarity. Here is nonetheless a summary.

We temporarily call $a$ the real part of $z$. The problem of finding $\widetilde{\Gamma}_+$ is divided into four parts. First, when $a \leq 1/\lambda_1$, we go along a line that makes an angle of $\pi/3$ with the real axis, starting at $c$. When $1/\lambda_1 \leq a \leq 1/\lambda_p$, we use a slightly more complicated path described in Section 3.4.2. When $1/\lambda_p$ is crossed, we go along a line that is parallel to the real axis until reaching a value $R_1$. At this value $R_1$, we go down vertically to the real axis.

Hence we will show that one can follow, even in the general case, a path that resembles that of Figure 4 in [6]. See Figure 2. There are two extra difficulties in the general case: we have to take care of an arbitrary spectrum which significantly increases the technical problems. Also, crossing the 1/eigenvalue zone is not a simple problem when first encountered.

In all that follows, we will use the notation

$$\alpha \triangleq c\lambda, \qquad \alpha_1 \triangleq c\lambda_1, \qquad \alpha_p \triangleq c\lambda_p, \qquad \gamma^2 \triangleq \frac{n}{p}, \qquad \mu \triangleq \mu_{n,p}.$$

We work under the assumptions of Theorem 1, hence $0 < c < 1/\lambda_1$, $0 < \alpha < 1$, $\limsup \alpha_1 < \infty$ and $\liminf \alpha_p > 0$. Recall that

$$f(z) = -\mu(z-q) + \log(z) - \frac{p}{n}\int \log(1-z\lambda)\, dH_p(\lambda).$$

3.4.1. *Behavior on* $\Gamma_1$. On $\Gamma_1$, we have $z = c + te^{i\pi/3}$.

We call $t = xc$ and consider $m(x) = \Re(f(c + xce^{i\pi/3}))$. Note that "$x$ increases" is equivalent to "$\Re(z)$ increases." This reparametrization considerably simplifies the computations. We have

$$\begin{aligned} m(x) = {} & -\mu c[(1+x/2) - q/c] + \frac{1}{2}\log(c^2(1+x+x^2)) \\ & - \frac{1}{2\gamma^2}\int \log((1-\alpha)^2 - x\alpha(1-\alpha) + \alpha^2 x^2)\, dH_p(\lambda). \end{aligned}$$

Recall that we want to show that $m'(x) < 0$, so that $m$ decreases when we move along $\Gamma_1$ with $\Re(z)$ (or equivalently $x$) increasing. We have

$$m'(x) = -\frac{\mu c}{2} + \frac{1}{2}\frac{2x+1}{1+x+x^2} - \frac{1}{2\gamma^2}\int \frac{2\alpha^2 x - \alpha(1-\alpha)}{(1-\alpha)^2 - x\alpha(1-\alpha) + \alpha^2 x^2}\, dH_p(\lambda).$$

Now remark that

$$\frac{n}{p} = \gamma^2 = \int \frac{\alpha^2}{(1-\alpha)^2}\, dH_p(\lambda)$$

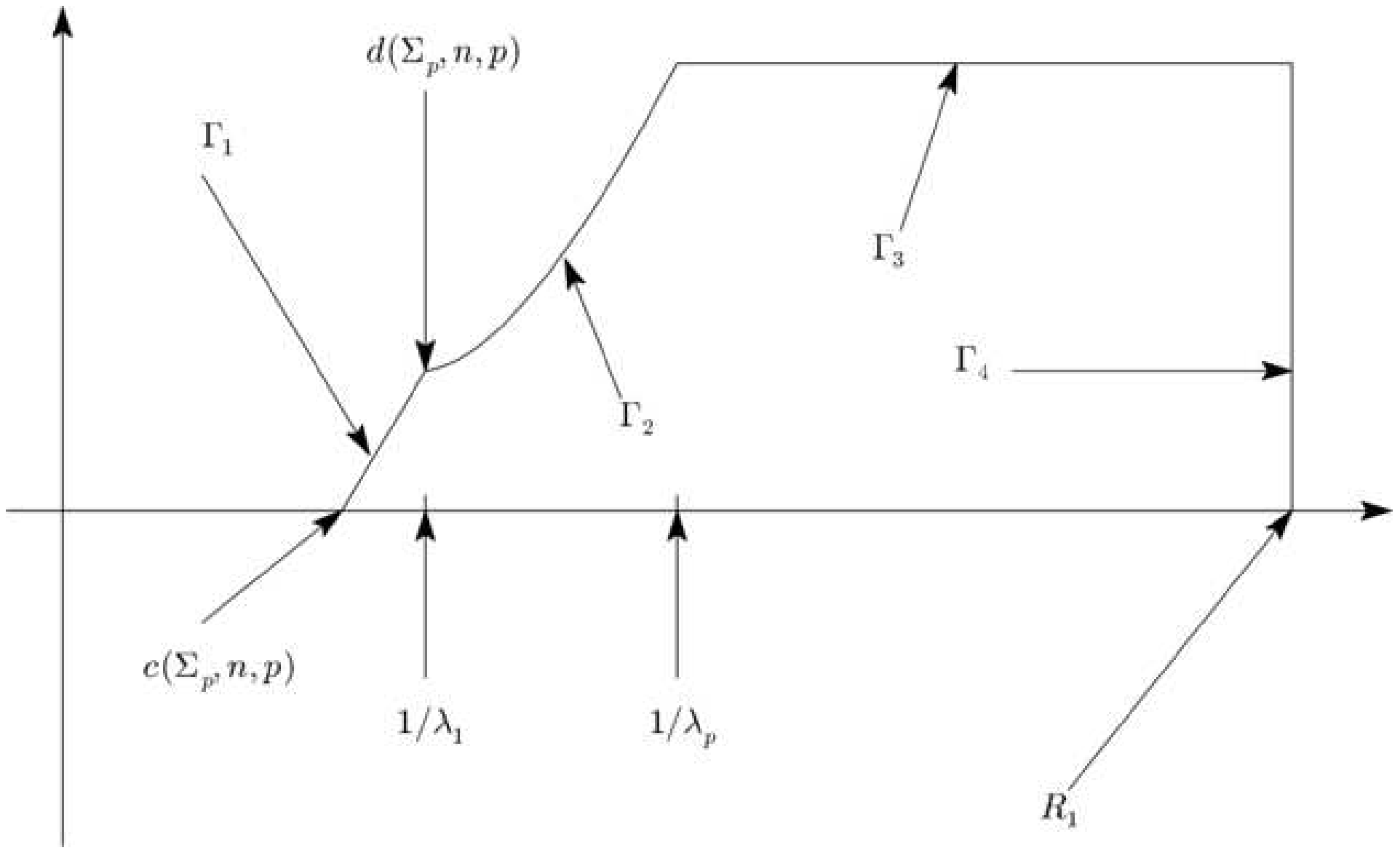

FIG. 2. *The curve* $\widetilde{\Gamma}_+$.

and

$$\begin{aligned}\mu c\gamma^2 &= \gamma^2 + \int \frac{\alpha}{1-\alpha}\, dH_p(\lambda) \\ &= \int \left( \frac{\alpha^2}{(1-\alpha)^2} + \frac{\alpha}{1-\alpha} \right) dH_p(\lambda) \\ &= \int \frac{\alpha^2 + \alpha(1-\alpha)}{(1-\alpha)^2}\, dH_p(\lambda) \\ &= \int \frac{\alpha}{(1-\alpha)^2}\, dH_p(\lambda).\end{aligned}$$

Therefore, $m'(x)2\gamma^2$ is equal to

$$\begin{aligned}&-\mu c\gamma^2 + \frac{2x+1}{1+x+x^2}\gamma^2 - \int \frac{2\alpha^2 x - \alpha(1-\alpha)}{(1-\alpha)^2 - x\alpha(1-\alpha) + \alpha^2 x^2}\, dH_p(\lambda) \\ &\quad = \int \biggl( -\frac{\alpha}{(1-\alpha)^2} + \frac{2x+1}{1+x+x^2} \frac{\alpha^2}{(1-\alpha)^2} \\ &\qquad\qquad - \frac{2\alpha^2 x - \alpha(1-\alpha)}{(1-\alpha)^2 - x\alpha(1-\alpha) + \alpha^2 x^2} \biggr) dH_p(\lambda) \\ &\quad = \int \frac{\alpha}{(1-\alpha)^2} \left[ -1 + \frac{\alpha(2x+1)}{1+x+x^2} - \frac{(1-\alpha)^2(2\alpha x - (1-\alpha))}{(1-\alpha)^2 - x\alpha(1-\alpha) + \alpha^2 x^2} \right] dH_p(\lambda).\end{aligned}$$

To simplify the problem, we note that the expression between the brackets can be written

$$\frac{\sum_{k=0}^{4} c_k x^k}{(1+x+x^2)((1-\alpha)^2 - x\alpha(1-\alpha)+\alpha^2x^2)}$$
$$= \frac{g(x,\alpha)}{(1+x+x^2)((1-\alpha)^2 - x\alpha(1-\alpha)+\alpha^2x^2)}.$$

A simple computation shows the following simplification:

$$c_0 = 0,$$
$$c_1 = 0,$$
$$c_2 = 2\alpha(\alpha-1),$$
$$c_3 = \alpha(2\alpha-1),$$
$$c_4 = -\alpha^2.$$

Hence

$$g(x,\alpha) = -\alpha x^2(2(1-\alpha)+(1-2\alpha)x+\alpha x^2).$$

We want $g(x,\alpha)$ to be negative, so we just have to study the polynomial $P(x,\alpha) = 2(1-\alpha)+(1-2\alpha)x+\alpha x^2$. Recall that $x \geq 0$. If $\alpha \leq 1/2$, all the coefficients are positive so the polynomial is positive for all $x \in \mathbb{R}_+$. Now the roots, at $\alpha$ fixed, of $P(\cdot,\alpha)$ are

$$x_{\pm} = \frac{(2\alpha-1) \pm \sqrt{(2\alpha-1)^2 - 8\alpha(1-\alpha)}}{2\alpha}.$$

The polynomial under the square root can be rewritten $h(\alpha) = 12\alpha^2 - 12\alpha + 1$. Its roots are $1/2 \pm 1/\sqrt{6}$, and it is negative between them.

Therefore, if $\alpha \leq 1/2 + 1/\sqrt{6} \simeq 0.9$, $P(x,\alpha) \geq 0$ for all $x$ in $\mathbb{R}_+$. So we just need to focus on $\alpha$'s such that $\alpha \geq 1/2 + 1/\sqrt{6}$.

Remark that $x_+$ and $x_-$ are both positive, because $\alpha \geq 1/2$. We potentially have a problem (of sign) when crossing the smaller one of the two roots, which is of course $x_-$. Now note that $x_-(\alpha) \geq x_-(\alpha_1)$.

As a matter of fact, we remark that $x_-(\alpha) < 1$ for all $\alpha$'s under consideration. Then, for $u \leq 1$, $P'(u,\alpha_1) \leq P'(u,\alpha)$, since $P'(u,\alpha) = 1 + 2\alpha(u-1)$ and $u \leq 1$. Note also that $P(0,\alpha_1) \leq P(0,\alpha) = 2(1-\alpha)$. So $P(u,\alpha_1) \leq P(u,\alpha)$ for $u \leq 1$. Since $P(x_-(\alpha_1),\alpha_1) = 0$, we see that $x_-(\alpha_1) \leq x_-(\alpha)$, if $\alpha \geq \alpha_1$. So for $u \leq x_-(\alpha_1)$, $P(u,\alpha) \geq 0$.

Now the only thing we need to verify to make sure that we can reach $\Re(z) = 1/\lambda_1$ is that $x_-(\alpha_1) \geq 2(\frac{1}{\alpha_1} - 1)$. Given that $\alpha_1 > 0.5 + 1/\sqrt{6}$, this is equivalent to showing that $\sqrt{(2\alpha_1 - 1)^2 - 8\alpha_1(1-\alpha_1)} \leq 6\alpha_1 - 5$, which is equivalent to $0 \leq 24(1-\alpha_1)^2$.

So we have shown that $\Re(f(z))$ decreases when the real part of $z$ increases, when going along the line intersecting the real axis at $c$ and making an angle of $\pi/3$ with it. If $\alpha_1 \leq 0.5 + 1/\sqrt{6}$, we can cross the whole plane along this line and $\Re(f(z))$ continues to decrease. If $\alpha_1 > 0.5 + 1/\sqrt{6}$, we are guaranteed that the property holds until $\Re(z) = 1/\lambda_1$.

Hence the claim we made about $\Gamma_1$ being a descent path of $\Re(f)$ is verified.

3.4.2. *Behavior on* $\Gamma_2$. As we saw in the previous subsection, this is only a concern if $\alpha_1 > 0.5 + 1/\sqrt{6}$. So we suppose we are in this situation. Before we proceed to exhibiting a path, we perform a preliminary computation that will prove useful in both this subsection and the next one.

*Independent computation.* Suppose we write $z = c(u + iv)$. We have

$$\Re(f(z)) = -\mu c(u - q/c) + \frac{1}{2}\log(c^2(u^2 + v^2)) - \frac{1}{2\gamma^2}\int \log((1 - \alpha u)^2 + \alpha^2 v^2)\, dH_p(\lambda).$$

If we consider that $v = v(u)$, we have $\Re(f(z)) = g(u)$. The question of finding a path along which $\Re(f(z))$ decreases when $\Re(z)$ increases is equivalent to finding $v(u)$ such that $g'(u) < 0$. With this in mind, we observe that

$$g'(u) = -\mu c + \frac{1}{2}\frac{2u + 2vv'}{u^2 + v^2} - \frac{1}{2\gamma^2}\int \frac{2\alpha(\alpha u - 1) + 2\alpha^2 vv'}{(1 - \alpha u)^2 + \alpha^2 v^2}\, dH_p(\lambda).$$

Let us call $I(u) = u^2 + v^2$ and $\beta = u + vv'$. Using the fact that $\mu c\gamma^2 = \int \alpha/(1 - \alpha)^2\, dH_p(\lambda)$ and $\gamma^2 = \int \alpha^2/(1 - \alpha)^2\, dH_p(\lambda)$, we get

$$\gamma^2 g'(u) = \int \frac{\alpha}{(1 - \alpha)^2}\left[-1 + \alpha\frac{\beta}{I(u)} - \frac{(\alpha\beta - 1)(1 - \alpha)^2}{(1 - \alpha u)^2 + v^2}\right] dH_p(\lambda). \tag{13}$$

*Back to the topic of* $\Gamma_2$. When $1/\lambda_1 \leq \Re(z) \leq 1/\lambda_p$, we have $1/\alpha_1 \leq u \leq 1/\alpha_p$. In this part of the plane, we propose to choose $\beta = I(u)$. Then the expression inside the brackets in equation (13) becomes

$$\begin{aligned}
-1 + \alpha &- \frac{(\alpha I(u) - 1)(1 - \alpha)^2}{\alpha^2 I(u) - 2\alpha u + 1} \\
&= (\alpha - 1)\left[1 - (\alpha - 1)\frac{\alpha I(u) - 1}{\alpha^2 I(u) - 2\alpha u + 1}\right] \\
&= (\alpha - 1)\frac{\alpha^2 I(u) - 2\alpha u + 1 - \alpha^2 I(u) + \alpha + \alpha I(u) - 1}{\alpha^2 I(u) - 2\alpha u + 1} \\
&= (\alpha - 1)\alpha\frac{I(u) - 2u + 1}{\alpha^2 I(u) - 2\alpha u + 1} \\
&= (\alpha - 1)\alpha\frac{(u - 1)^2 + v^2}{(1 - \alpha u)^2 + v^2} < 0.
\end{aligned}$$

Note that at the end of $\Gamma_1$ we arrived at $u_1 = 1/\alpha_1$ and the corresponding $v$ was $v_1 = \frac{\sqrt{3}(1-\alpha_1)}{\alpha_1}$. Now the choice of $\beta = I(u)$ can be reformulated as $I'(u) = 2I(u)$ and hence $I(u) = Ke^{2u}$. Simple algebra shows that finally

$$I(u) = \left(\frac{1}{\alpha_1}\right)^2 (1 + 3(1-\alpha_1)^2)e^{2u-2/\alpha_1} \qquad \text{on } \Gamma_2.$$

Note also that since $u \geq 1$, $I(u) = u^2 + v^2 > u$ and hence $\beta = u + vv' = I(u)$ implies that $v' > 0$, as we started with $v_1 > 0$. So we will not cross the real axis by following this path. In the original coordinates, if we call $z = a + ib$, the path is such that

$$b^2 = \left(\frac{1}{\lambda_1}\right)^2 (1 + 3(1-\alpha_1)^2)e^{2(a-1/\lambda_1)/c} - a^2,$$

with $b > 0$. For $\Gamma_2$, we follow this path until we reach $a = 1/\lambda_p$. Note that with our assumptions about $\gamma$, $\limsup \alpha_1$ and $\liminf \alpha_p$, the length of this path is uniformly bounded. We also remark that if $\alpha_p \to 0$, the length of $\Gamma_2$ grows to $\infty$, which causes problem for the control of the operator later on.

3.4.3. *Behavior on* $\Gamma_3$. We revert to the notation $z = c(u + iv)$. The point is just to show that with $v' = 0$ when $u > 1/\alpha_p$, $\Re(f(z))$ decreases. If we recall (13), we realize that if $\alpha\beta \leq I(u)$ and $\alpha\beta \geq 1$, then $g'(u) \leq 0$. But when $v' = 0$, $\beta = u$. Now, if $u \geq 1/\alpha_p$, $\alpha u = \alpha\beta \geq \alpha_p\beta \geq 1$. Also, since $u \leq I(u)$, and $\alpha \leq 1$, $\alpha u \leq I(u)$. Hence $\Re(f(z))$ is decreasing when moving along $\Gamma_3$.

3.4.4. *Behavior of* $\Gamma_4$. There, $z = R_1 + iy$, where, with a slight abuse of notation, $0 < y < \Gamma_2(1/\alpha_p)$. Now

$$\begin{aligned}\Re(f(z)) &= -\mu(R_1 - q) + \frac{1}{2}\log(R_1^2 + y^2) \\ &\qquad - \frac{1}{2\gamma^2}\int \log((\lambda R_1 - 1)^2 + \lambda^2 y^2)\, dH_p(\lambda)\end{aligned}$$

so, since $\mu$ is bounded away from 0, if $R_1 \to \infty$, $\Re(f(z)) \to -\infty$, and we can pick $R_1$ so that, uniformly for our models,

$$\Re(f(z)) \leq \Re(f(d)).$$

This is a simple consequence of the fact that $\Re(f(d(\Sigma_p, n, p)))$ is bounded below under the assumptions of Theorem 1. (See Appendix A.1.1.)

There is of course a problem of definition of $f$ at $y = 0$, because the argument of the logarithm is real and negative. Nevertheless the function $h(z) = e^{-n\mu_{n,p}(z-q)} z^n \prod_{k=1}^p (1 - z\lambda_k)^{-1}$ is well defined and well behaved at $z = R_1$, so this definition problem will cause no harm in the analysis of the convergence of the operators. As a matter of fact, it turns out that we will just be interested in bounding $|h(z)|$. Since we can take the log of $|h(z)|$ without any problems and it leads to the same expression as the one for $\Re(f(z))$ we considered, we can safely ignore the definition of $f$ problem for all practical purposes.

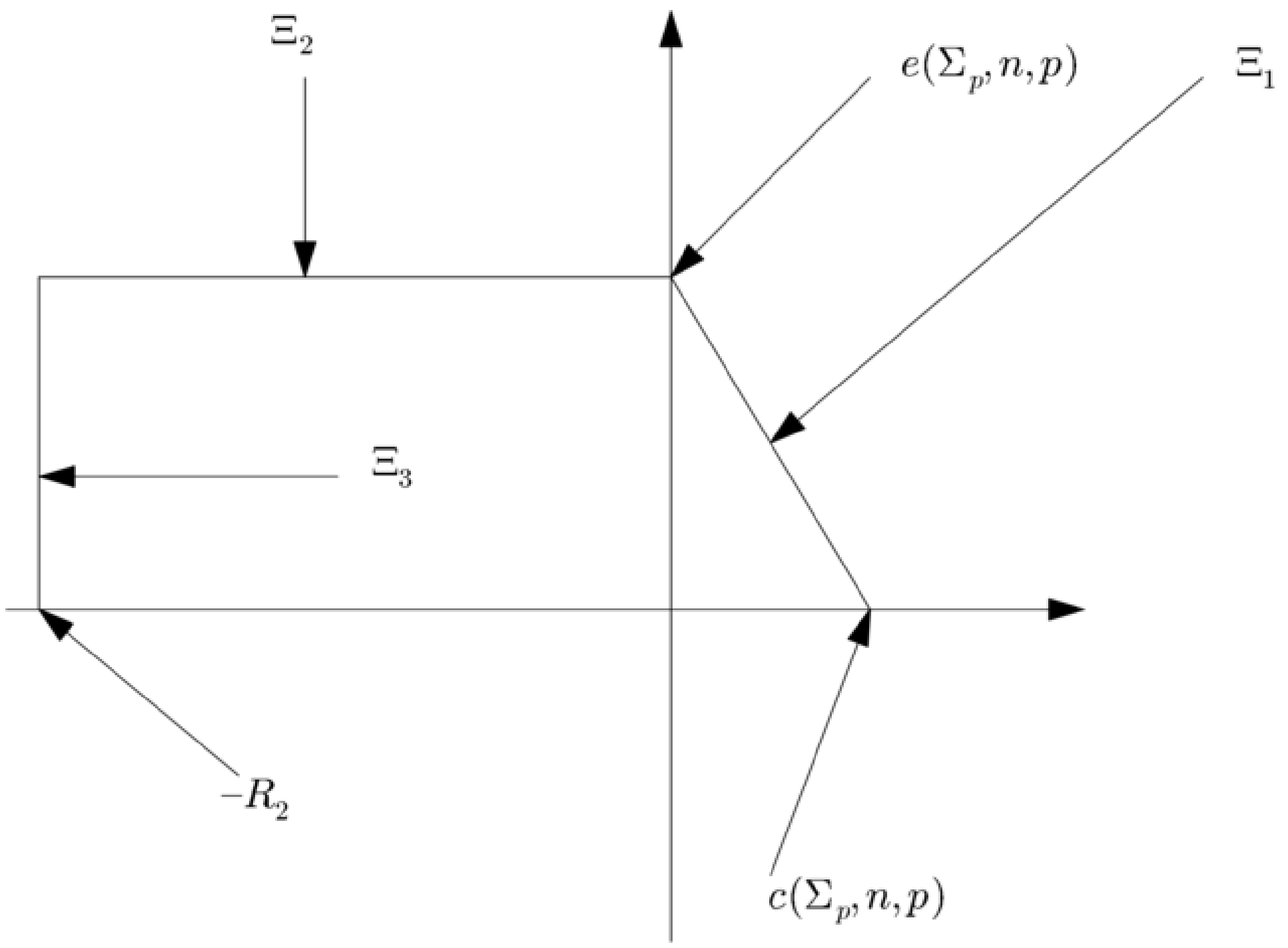

FIG. 3. *The curve* $\tilde{\tilde{\Xi}}_+$.

3.5. *About* $\Xi$. We use the same conventions as when we studied $\Gamma$. Namely, we will study the behavior of $f$ on $\Xi$, but we will first exhibit $\tilde{\tilde{\Xi}}$, with $\tilde{\tilde{\Xi}} = \tilde{\tilde{\Xi}}_+ \cup \overline{\tilde{\tilde{\Xi}}_+}$. It turns out that the analysis is much simpler for this contour and we will be able to follow the path used in [6], after doing some precise technical work.

Once again, the problem is very graphical. A drawing of $\tilde{\tilde{\Xi}}_+$ is shown in Figure 3.

What we will have to do in this case is to show that $\Re(-f(z))$ is decreasing when we travel along $\Xi_1$ and $\Xi_2$, and $\Re(z)$ is decreasing.

This time $\Xi_1$ is defined as a line making an angle of $2\pi/3$ with the real axis and crossing it at $c$. $\Xi_2$ is a line that runs parallel to the real axis, in the direction of $-\infty$.

The aim of this subsection is to show the following lemma:

LEMMA 2. *Under the assumptions of Theorem 1, $\Re(-f(z))$ is decreasing for $z \in \Xi_1 \cup \Xi_2$ as $\Re(z)$ decreases. Also, the length of $\widetilde{\Xi}_+$ is uniformly bounded. Finally, there exists $R_2 > 0$ such that $\max_{z\in\Xi_3} \Re(-f(z)) \le \Re(-f(e))$, where $e = e(\Sigma_p, n, p) = ic\sqrt{3}$.*

3.5.1. *Case of* $\Xi_1$. Once again, we will consider everything on the $c$ scale. We define $\Xi_1$ as $z = c + xce^{i2\pi/3}$. We have

$$\begin{aligned}
\phi_1(x) &\triangleq \Re(-f(c(1 + xe^{i2\pi/3}))) \\
&= \Re(-f(c(1 - x/2 + ix\sqrt{3}/2))) \\
&= \mu c(1 - x/2) - \frac{1}{2}\log(1 - x + x^2) \\
&\quad + \frac{1}{2\gamma^2}\int \log((1 - \alpha(1 - x/2))^2 + 3\alpha^2x^2/4)\, dH_p(\lambda) \\
&\quad + \frac{1}{2}\left(\frac{1}{\gamma^2} - 1\right)\log(c^2) + \mu q.
\end{aligned}$$

Hence, we get

$$\begin{aligned}
\gamma^2\phi_1'(x) &= -\frac{1}{2}\mu c\gamma^2 - \gamma^2\frac{2x - 1}{2(1 - x + x^2)} \\
&\quad + \frac{1}{2}\int \frac{2\alpha^2 x + \alpha(1 - \alpha)}{(1 - \alpha)^2 + x\alpha(1 - \alpha) + \alpha^2x^2}\, dH_p(\lambda).
\end{aligned}$$

Therefore, using the same equalities we used when studying $\Gamma$, we have

$$\begin{aligned}
2\gamma^2\phi_1'(x) &= \int \frac{-\alpha}{(1 - \alpha)^2} - \frac{2x - 1}{1 - x + x^2}\frac{\alpha^2}{(1 - \alpha)^2} \\
&\qquad + \frac{2\alpha^2x + \alpha(1 - \alpha)}{\alpha^2x^2 + \alpha(1 - \alpha)x + (1 - \alpha)^2}\, dH_p(\lambda) \\
&= \int \frac{\alpha}{(1 - \alpha)^2}\biggl(-1 - \alpha\frac{2x - 1}{1 - x + x^2} \\
&\qquad\qquad + \frac{(2\alpha x + (1 - \alpha))(1 - \alpha)^2}{\alpha^2x^2 + \alpha(1 - \alpha)x + (1 - \alpha)^2}\biggr)\, dH_p(\lambda).
\end{aligned}$$

As before, the expression that is within the parentheses can be written

$$\frac{\sum_{k=0}^4 c_k x^k}{((1 - x/2)^2 + 3x^2/4)(\alpha^2x^2 + \alpha(1 - \alpha)x + (1 - \alpha)^2)}$$

and we know that the denominator is positive. A simple computation leads to

$$c_0 = 0,$$

$$\begin{aligned} c_1 &= 0, \\ c_2 &= -2\alpha + 2\alpha^2, \\ c_3 &= \alpha - 2\alpha^2, \\ c_4 &= -\alpha^2, \end{aligned}$$

and hence the numerator is

$$\sum_{k=0}^{4} c_k x^k = -x^2\alpha(\alpha x^2 + (2\alpha - 1)x + 2(1-\alpha)).$$

The same questions we asked when dealing with $\Gamma_1$ now come up. Note that our $x$ is positive, so if $\alpha \geq 1/2$, $P(x,\alpha) = (\alpha x^2 + (2\alpha - 1)x + 2(1-\alpha)) \geq 0$. Also, at $\alpha$ fixed the roots of $P(\cdot,\alpha)$ are

$$x_{\pm} = \frac{1 - 2\alpha \pm \sqrt{12\alpha^2 - 12\alpha + 1}}{2\alpha}.$$

As we saw before, we therefore have $P(x,\alpha) \geq 0$ on $\mathbb{R}_+ \times [1/2 - 1/\sqrt{6}, 1]$. Now if $\alpha \leq 1/2 - 1/\sqrt{6}$, we have to work a little harder. We remark that if $\alpha \in [0, 0.5 - 1/\sqrt{6}]$, it is easy to check that $x_-(\alpha) \geq 2$. Hence we conclude that

$$P(x,\alpha) \geq 0 \qquad \text{on } [0,2] \times [0,1].$$

Now $z(2) = 0 + ic\sqrt{3} = e$. So we have shown that $\Re(-f(z))$ decreases when we travel from $c$ to $e$ along $\Xi_1$.

3.5.2. *Case of* $\Xi_2$. On this part of the path we use the parametrization $z = -xc + i\sqrt{3}c$, with $x \geq 0$ and increasing. We have, if $K$ is a constant (at $\Sigma_p$, $n$ and $p$ fixed),

$$\Re(f(z)) = \mu xc + \frac{1}{2}\log(x^2+3) - \frac{1}{2\gamma^2}\int \log((1+\alpha x)^2 + 3\alpha^2)\, dH_p(\lambda) + K.$$

Calling $\phi_2(x) = \Re(-f(z))$, we hence get

$$\phi_2(x) = -\mu xc - \frac{1}{2}\log(x^2+3) + \frac{1}{2\gamma^2}\int \log((1+\alpha x)^2 + 3\alpha^2)\, dH_p(\lambda) + K.$$

Using the same approach as before we find that

$$\gamma^2\phi_2'(x) = \int \frac{\alpha}{(1-\alpha)^2}\left[-1 - \alpha\frac{x}{x^2+3} + \frac{(1+\alpha x)(1-\alpha^2)}{(1+\alpha x)^2 + 3\alpha^2}\right] dH_p(\lambda).$$

Once again what matters to us is the numerator of what is within the bracket. It is a polynomial—let us denote it by $Q(x,\alpha)$—of degree 4 in

$x$, its coefficients being

$$\begin{aligned} c_0 &= -6\alpha(1-\alpha), \\ c_1 &= -2\alpha(2+3\alpha), \\ c_2 &= -\alpha(2+7\alpha), \\ c_3 &= -\alpha(1+2\alpha), \\ c_4 &= -\alpha^2. \end{aligned}$$

Hence it is clear that $Q(x,\alpha) \leq 0$ on $\mathbb{R}_+ \times [0,1]$. Therefore, we have shown that $\Re(-f(z))$ decreases as $\Re(z)$ decreases and $z$ travels on $\Xi_2$.

3.5.3. *Case of* $\Xi_3$. Here $z = -R_2 + iy$, where $0 \leq y \leq c\sqrt{3}$. It is easy to see that $\Re(-f(z))$ can be made as small as we want, since $\mu$ is bounded away from 0. We show in Appendix A.1.2 that $\Re(-f(e))$ is bounded. In particular, if we choose $R_2$ large enough,

$$\max_{z\in\Xi_3} \Re(-f(z)) \leq \Re(-f(e)).$$

This holds uniformly with respect to our covariance models, if they are in $\mathcal{G}$.

3.6. *Study of* $A_{n,p}$. We give an outline of the key ideas and results that allow us to then proceed to operator convergence issues. The proof is given in Appendix B.

3.6.1. *Definition of* $q$ *and modification of* $\widetilde{\Gamma}$ *to get* $\Gamma$. At this point we still have not set $q$ and it is now time to do it. Let us pick an $\varepsilon > 0$. Then, set

$$q \triangleq q(\Sigma_p, n, p) = c - \frac{\varepsilon}{n\sigma_{n,p}}. \tag{14}$$

Then, as in [6], we just have to modify the curve $\widetilde{\Gamma}_+$ around $c$ to obtain $\Gamma_+$. $\Gamma_+$ is the same as $\widetilde{\Gamma}_+$, except it starts by

$$\Gamma_0 = \left\{ c + \frac{\varepsilon}{2n\sigma_{n,p}} e^{i\theta} : \frac{\pi}{3} \leq \theta \leq \pi \right\}.$$

When $\Gamma_0$ reaches $\Gamma_1$, we follow $\Gamma_1$, and then follow $\Gamma_2$, $\Gamma_3$ and $\Gamma_4$ to create $\Gamma_+$. Then $\Gamma = \Gamma_+ \cup \overline{\Gamma_+}$. Of course, in the end, the contour $\Gamma$ is oriented counterclockwise. A depiction of $\Gamma_+$ can be found in Figure 4.

3.6.2. *Arguments needed for the operator analysis to go through.* The method of proof is similar to that of [6], once the difficulties stemming from the fact that we are considering a much more general case are understood.

The issue we will face is to find a sequence $\kappa_{n,p}$ such that $\kappa_{n,p}A_{n,p} \to e^{-\varepsilon x}\operatorname{Ai}(x)$ and $\kappa_{n,p}A_{n,p}$ goes to zero exponentially fast at infinity.

The analysis will rely on four key points. They are:

1. The length of $\Gamma$ is uniformly bounded with respect to our models. We will justify in Appendix B.1 why this is the case in the situation we are considering.
2. One needs to be able to find $\delta > 0$ such that
$$\forall s: |s-c| < \delta \quad \Longrightarrow \quad \frac{|f^{(4)}(s)|}{4!}\delta < \frac{\sigma^3}{6}.$$
$\delta$ has of course to be uniform with respect to our models.
3. We also need
$$\forall s: |s-c| < \delta \quad \Longrightarrow \quad \limsup \sup \frac{|f^{(4)}(s)|}{4!} = \Delta < \infty.$$
4. Finally, $\delta$ has to be chosen small enough that the disc of center $c$ and radius $\delta$ should encompass neither $d(\Sigma_p, n, p)$ nor $e(\Sigma_p, n, p)$.

We will explain in Appendix B.1 why these conditions are fulfilled under the assumptions of Theorem 1 and then prove in Appendix B.2 the following proposition:

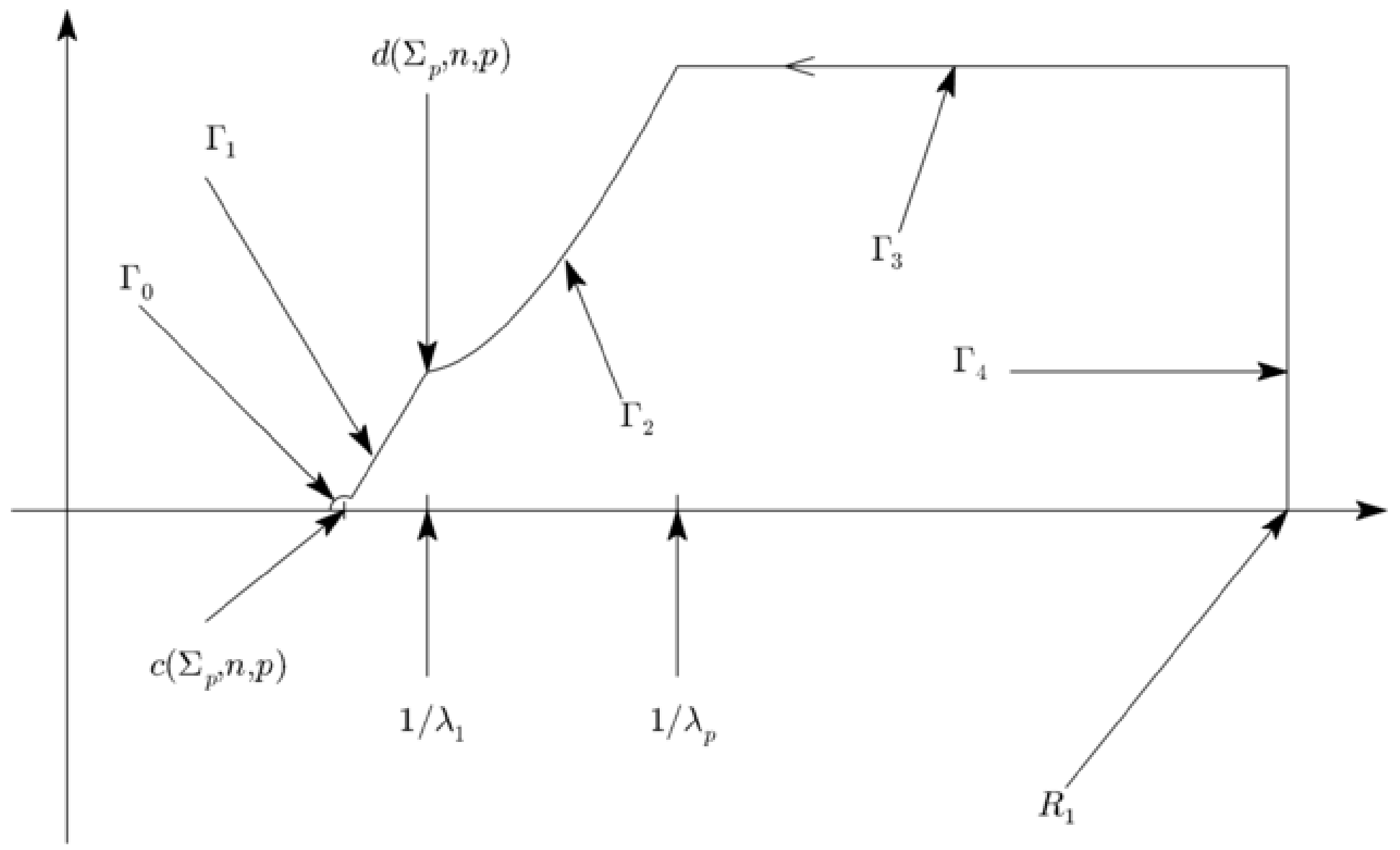

FIG. 4. *The curve $\Gamma_+$.*

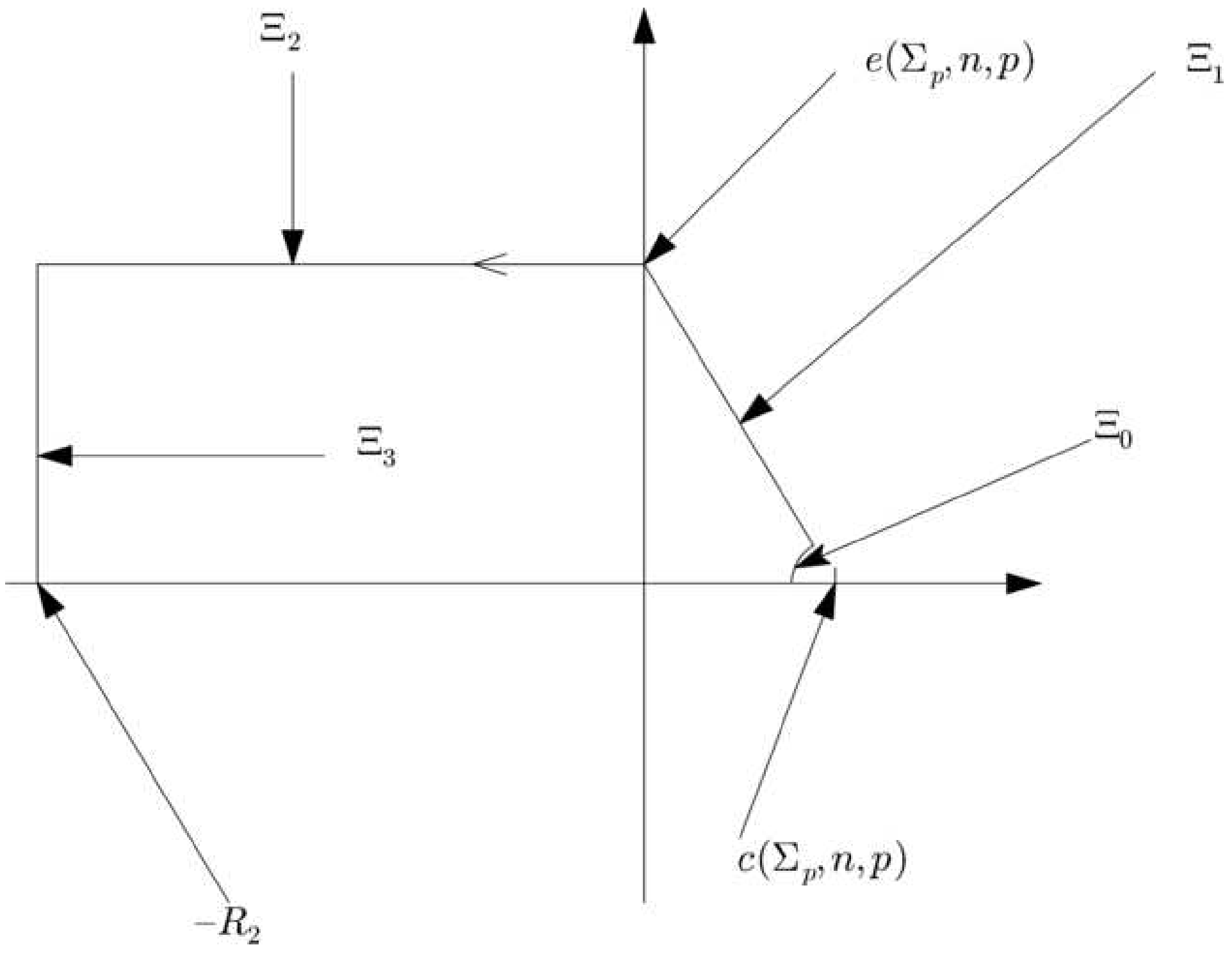

FIG. 5. *The curve* $\Xi_+$.

PROPOSITION 1. *In the definition* [*see* (9)] *of* $A_{n,p}$, *let* $\mu_{n,p}$ *be equal to* $\mu$ *in* (2) *and* $\sigma_{n,p} = n^{-2/3}\sigma$, *with* $\sigma$ *defined in* (3). *When the four conditions above are fulfilled, we have*

$$\forall s_0 \in \mathbb{R}, \exists C(s_0) \in \mathbb{R}_+ \text{ and } N_0 \in \mathbb{N} \text{ such that}$$

$$|\kappa_{n,p} A_{n,p}(s) - e^{-\varepsilon s}\mathrm{Ai}(s)| \leq \frac{C(s_0)e^{-\varepsilon s/2}}{n^{1/3}}$$

*if* $s \geq s_0$ *and* $n \geq N_0$. *Here* $\kappa_{n,p} = e^{-nf(c)}$.

*As a function of* $s_0$, $C$ *can be chosen to be continuous and nonincreasing.*

3.7. *Study of* $B_{n,p}$. Here also, $\widetilde{\Xi}$ needs to be modified. We start by $\Xi_0$, an arc of a circle centered at $c$ and with radius $3\varepsilon/(n\sigma_{n,p})$. Formally, $\Xi_0 = c + 3\varepsilon/(n\sigma_{n,p})e^{i(\pi-\theta)}$, with $0 \leq \theta \leq \pi/3$. When $\Xi_0$ intersects $\Xi_1$, we follow $\Xi_1$, and so on. A depiction of $\Xi_+$ can be found in Figure 5.

We then have:

PROPOSITION 2. *In the definition* [*see* (10)] *of* $B_{n,p}$, *let* $\mu_{n,p}$ *be equal to* $\mu$ *in* (2) *and* $\sigma_{n,p} = n^{-2/3}\sigma$, *with* $\sigma$ *defined in* (3). *When the four conditions* *in Section 3.6.2 are fulfilled,*

$$\forall s_0 \in \mathbb{R}, \exists C(s_0) \in \mathbb{R}_+ \text{ and } N_0 \in \mathbb{N} \text{ such that}$$

$$|B_{n,p}(s)/\kappa_{n,p} - e^{\varepsilon s}\mathrm{Ai}(s)| \leq \frac{C(s_0)e^{-\varepsilon s/2}}{n^{1/3}}$$

*if* $s \geq s_0$ *and* $n \geq N_0$. *Here* $\kappa_{n,p} = e^{-nf(c)}$, *the same as in Proposition 1.*

*As a function of* $s_0$, $C$ *can be chosen to be continuous and nonincreasing.*

Explanations are postponed to Appendix B.3.

3.8. *Operator convergence issues.* We will mostly rely on two key properties in this subsection: the relationship between trace class and Hilbert–Schmidt norms and the fact that the determinant $\det(I - \cdot)$ of trace class operators is a locally Lipschitz function with respect to trace class norm.

More precisely, recall (see [16], Section IV.7) that if $O$ and $P$ are Hilbert–Schmidt operators, then $OP$ is a trace class operator and

$$\|OP\|_1 \leq \|O\|_2\|P\|_2.$$

Also, it is well known (see [16], Theorem IV.5.2, and Theorem II.4.1 and Corollary II.4.2 both due to Seiler–Simon) that if $Q$ and $R$ are trace class operators,

$$\text{(Lip)} \qquad |\det(I+Q) - \det(I+R)| \leq \|Q-R\|_1 e^{\|Q\|_1 + \|R\|_1 + 1}.$$

This section is now devoted to proving two lemmas that allow us to prove Theorem 1.

Let us call $E$ the multiplication operator by $e^{-x}$ and $\mathrm{Ai}_s$ the operator on $L^2([s,\infty))$ with kernel $\mathrm{Ai}(x+y-s)$.

Lemma 3. *Using the conclusions of Propositions 1 and 2, we have, if we view all the operators as operators on* $L^2([s,\infty))$*:*

$$\forall s_0 \in \mathbb{R}, \exists B \in \mathbb{R}_+ \text{ and } N_0 \in \mathbb{N} \text{ such that}$$

$$\|A_{n,p}B_{n,p} - E\mathrm{Ai}_s^2 E^{-1}\|_1 \leq \frac{C(s_0)e^{-\varepsilon s/2}}{n^{1/3}},$$

*if* $s \geq s_0$ *and* $n \geq N_0$. $C$, *as a function of* $s_0$, *can be chosen to be continuous and nonincreasing.*

Proof. Recall the following fact: according to [28], page 394, for $x > 0$, $\mathrm{Ai}(x) \leq \exp(-2x^{3/2}/3)/(2\pi^{1/2}x^{1/4})$. Hence it is clear that the operator $P$ with kernel $P(x,y) = P(x+y-s) = \mathrm{Ai}(x+y-s)\exp(\varepsilon(x+y-s))$ is Hilbert–Schmidt on $L^2([s,\infty))$, and similarly for $O$ that has kernel $O(x,y) = O(x+y-s) = \mathrm{Ai}(x+y-s)\exp(-\varepsilon(x+y-s))$.

More precisely, since these kernels are as functions of $(x,y)$ square integrable on $[s,\infty)\times[s,\infty)$, Theorem VI.23 in [30] applies and we see that, for instance, if we view $O$ as an operator on $L^2([s,\infty))$,

$$\begin{aligned}\|O\|_2^2 &= \iint_{[s,\infty)^2} (O(x+y-s))^2\,dx\,dy\\ &= \iint_{[0,\infty)^2} (O(x+y+s))^2\,dx\,dy\\ &= \int_{x=s}^{\infty}\int_{y=0}^{\infty} (O(x+y))^2\,dx\,dy.\end{aligned}$$

It is clear that this is a continuous, nonincreasing function of $s$ having limit 0 at $\infty$. The same analysis and conclusion apply to $P$.

Now let us denote $\widetilde{A}_{n,p}=\kappa_{n,p}A_{n,p}$ and $\widetilde{B}_{n,p}=B_{n,p}/\kappa_{n,p}$. From the previous analyses we conclude that we can find a continuous, nonincreasing function $C$ such that, if we view all the operators as operators on $L^2([s,\infty))$, with $s\ge s_0$, $\|\widetilde{A}_{n,p}\|_2\le C(s_0)$, $\|P\|_2\le C(s_0)$, and similarly for $\widetilde{B}_{n,p}$ and $O$, as long as $n\ge N_0(s_0)$. For instance, $C(s)$ could be $2(\|O\|_2(s)+\|P\|_2(s))$, where we have highlighted the dependence of the Hilbert–Schmidt norm of $O$ and $P$ on $s$.

Since $A_{n,p}B_{n,p}-OP=\widetilde{A}_{n,p}(\widetilde{B}_{n,p}-P)+(\widetilde{A}_{n,p}-O)P$, we have

$$\|A_{n,p}B_{n,p}-OP\|_1\le\|\widetilde{A}_{n,p}\|_2\|\widetilde{B}_{n,p}-P\|_2+\|\widetilde{A}_{n,p}-O\|_2\|P\|_2.$$

Using the estimates obtained in Propositions 1 and 2, we have shown that if we view $A_{n,p}B_{n,p}-OP$ as an operator on $L^2([s,\infty))$ with $s\ge s_0$, we have

$$\|A_{n,p}B_{n,p}-OP\|_1\le\frac{C(s_0)\exp(-\varepsilon s/2)}{n^{1/3}},$$

if $n\ge N_0(s_0)$, for yet another continuous, nonincreasing function $C$. Finally, $OP$ and $E\mathrm{Ai}_s^2E^{-1}$ have the same kernel, so Lemma 3 is proved. □

Let us call $F_0$ the cumulative distribution function of the Tracy–Widom distribution arising in the study of the Gaussian unitary ensemble. Recall that, as explained in [35], formula (4.5) and page 166, $F_0(s)=\det(I-\mathrm{Ai}_s^2)$, where $\mathrm{Ai}_s^2$ is viewed as an operator on $L^2([s,\infty))$. Note that since $E\mathrm{Ai}_s$ and $\mathrm{Ai}_sE^{-1}$ are clearly Hilbert–Schmidt on $L^2([s,\infty))$, $\det(I-E\mathrm{Ai}_s\mathrm{Ai}_sE^{-1})=\det(I-\mathrm{Ai}_sE^{-1}E\mathrm{Ai}_s)=\det(I-\mathrm{Ai}_s^2)$.

We also have

$$P\biggl(\frac{l_1-n\mu}{\sigma n^{1/3}}\biggr)=\det(I-A_{n,p}B_{n,p}).$$

The continuity of the determinant $\det(I-\cdot)$ with respect to trace class norm implies that

$$\biggl|P\biggl(\frac{l_1-n\mu}{n^{1/3}\sigma}\le s\biggr)-F_0(s)\biggr|\to 0.$$

The convergence part of Theorem 1 is therefore proved.

We now turn to proving the rate of convergence part of Theorem 1. In other words, we want to show that:

*We can find a function* $C$ (*continuous and nonincreasing if we wish*) *such that* $\forall s_0, \exists N_0$:

$$s \geq s_0 \text{ and } n \geq N_0 \text{ implies}$$
$$\left| P\left( \frac{l_1 - n\mu}{n^{1/3}\sigma} \leq s \right) - F_0(s) \right| \leq \frac{C(s_0)e^{-\varepsilon s/2}}{n^{1/3}}.$$

PROOF. First, it is clear that since, in the notation of the previous proof, $\widetilde{A}_{n,p}$ and $\widetilde{B}_{n,p}$ are Hilbert–Schmidt operators and converge to, respectively, $O$ and $P$, we have, when considering our operators as operators on $L^2([s, \infty))$,

$$\forall s \geq s_0 \text{ and } n \text{ large enough } \|A_{n,p}B_{n,p}\|_1 \leq \|\widetilde{A}_{n,p}\|_2 \|\widetilde{B}_{n,p}\|_2 \leq 2\|O\|_2\|P\|_2.$$

This last quantity is less than $C(s)$, where $C$ is a continuous, nonincreasing function, going to 0 when $s$ tends to $\infty$.

Hence, for $s \geq s_0$, if $n > N_0(s_0)$, $\|A_{n,p}B_{n,p}\|_1 + \|OP\|_1 \leq 3C(s_0)$, for yet another continuous, nonincreasing function $C$. In view of equation (Lip) and the estimate we already have for $\|A_{n,p}B_{n,p} - OP\|_1$, the statement is shown, because $|P((l_1 - n\mu)/(n^{1/3}\sigma) \leq s) - F_0(s)| = |\det(I - A_{n,p}B_{n,p}) - \det(I - OP)|$.

Since the $C$'s appearing in Propositions 1 and 2 may depend on the models under consideration, so may $C$. □

Hence the rate of convergence part of Theorem 1 is proved.

**4. Simulations, related issues and conclusion.** We will discuss in this section some practical consequences of Theorem 1 as well as some of the questions it raises. To simplify the discussion, we recall that we denote by $\mathcal{G}$ the class of (covariance) models for which Theorem 1 applies. We will often abuse the notation and say that a covariance matrix is in $\mathcal{G}$ to mean that the corresponding model is in $\mathcal{G}$.

4.1. *Finite perturbation of a covariance matrix that is in* $\mathcal{G}$. In this subsection, we discuss some immediate consequences of the analysis we made to the case of a finite perturbation of a covariance matrix that is in $\mathcal{G}$. By this we mean that we are now considering data matrices $X$ that are $n \times (p + k)$, and $X_i \overset{\text{i.i.d.}}{\sim} \mathcal{N}_{\mathbb{C}}(0, \widetilde{\Sigma}_{p+k})$, where $k(p) < K$, $K \in \mathbb{N}$, and we add to $\{\lambda_1(\Sigma_p), \ldots, \lambda_p(\Sigma_p)\}$ $k$ eigenvalues larger than $\lambda_1(\Sigma_p)$. In other words, $\lambda_{k+1}(\widetilde{\Sigma}_{p+k}) = \lambda_1(\Sigma_p)$. This is of course a generalization of spiked covariance models considered in [6], where the bulk covariance is not restricted to be Id but rather a matrix for which Theorem 1 applies.

We will discuss two cases. First, we will assume that there exists $\chi > 0$ such that $\lambda_1(\widetilde{\Sigma}_{p+k}) < -\chi + 1/c(\Sigma_p, n, p)$. We will see in that case that Theorem 1 applies. Then we will discuss the case where $\lambda_1(\widetilde{\Sigma}_{p+k}) = 1/c(\Sigma_p, n, p)$ and the situation when the multiplicity of this eigenvalue is $k_0$.

FACT 1. *In the spiked situation described above, if there exists $\chi > 0$ such that*

$$\lambda_1(\widetilde{\Sigma}_{p+k}) < -\chi + 1/c(\Sigma_p, n, p),$$

*Theorem 1 applies to $\{(\widetilde{\Sigma}_{p+k}, n, p+k)\}$.*

The proof is elementary and is given in Appendix A.4.1. Intuitively this means that if we perturb a model for which Theorem 1 applies by adding a few leading eigenvalues that are not too large [and too large means larger than $1/c(\Sigma_p, n, p) - \chi$ for some $\chi > 0$], then Theorem 1 applies to the perturbed model.

In light of [6], another natural question is to understand what happens when we spike the model by adding $k$ eigenvalues at exactly $1/c(\Sigma_p, n, p)$. We have the following result in this case:

THEOREM 3. *Let us assume that a model in $\mathcal{G}$ is spiked by adding $k$ eigenvalues at*

$$\lambda_1(\widetilde{\Sigma}_{p+k}) = \cdots = \lambda_k(\widetilde{\Sigma}_{p+k}) = 1/c(\Sigma_p, n, p).$$

*The value of $k$ is fixed and is not allowed to change with $n$ or $p$. Then calling $F_i$'s the distribution functions defined in Definition 1.1 of [6], we have*

$$P\left(\frac{l_1 - n\mu}{n^{1/3}\sigma} \le x\right) \to F_k(x).$$

*As in Theorem 1, we have*

$$\mu = \frac{1}{c}\left(1 + \frac{p}{n}\int \frac{\lambda c}{1 - \lambda c}\, dH_p(\lambda)\right),$$

$$\sigma^3 = \frac{1}{c^3}\left(1 + \frac{p}{n}\int \left(\frac{\lambda c}{1 - \lambda c}\right)^3 dH_p(\lambda)\right).$$

*Note that $c$, $\mu$ and $\sigma$ refer to the nonspiked model.*

A justification is given in Appendix B. So we have extended Theorem 1.1(a) in [6] to models in the class $\mathcal{G}$. More information about the $F_k$'s can be found in [5], [6] and Appendix B.1.4.

4.2. *Statistical considerations.*

4.2.1. *Isolated largest eigenvalue vs. largest eigenvalue with a small mass.* One of the many very interesting results obtained in [6] was their Theorem 1.1(b). It basically says that if an Id matrix is spiked with eigenvalues that are larger than $1/c(\Sigma_p, n, p) + \chi$, $\chi > 0$, $l_1$ has a completely different type of limiting distribution, and that centering and scaling should be changed. In particular the scaling should be adjusted from $n^{1/3}$ to $n^{1/2}$. The question of knowing if and how this happens for matrices of the class $\mathcal{G}$ is currently under investigation by the author of this article. As an aside, let us remark that $n^{1/2}$ is the rate obtained through elementary concentration of measure arguments. We refer to Appendix A.5 and references therein for more details.

Let us go back to our discussion and call this large spike $\tilde{\lambda}_1$. If instead of changing one eigenvalue we had a small mass $\nu(p)$ [with $\liminf \nu(p) > 0$] at $\tilde{\lambda}_1$, then Theorem 1 would apply. Hence the centering, scaling and limiting distribution of $l_1$ would differ drastically from the case where $\tilde{\lambda}_1$ is isolated. In practice (and in statistical applications), one cannot tell from the data if there is one eigenvalue (out of say 100) that is much larger than the rest of them, or if 1% of the eigenvalues are clearly separated from the bulk. One will therefore have to specify precisely what models are considered if the results presented in this paper and those in [6] are used for statistical inference. Note that asymptotics done at fixed spectral distribution lead to Tracy–Widom limits.

For instance, in a hypothesis testing context, the power of tests based on these "large $p$, large $n$" asymptotics will depend greatly on the specified alternatives.

4.2.2. *Classical asymptotics or* $\lim n/p < \infty$ *asymptotics*? An interesting statistical aspect of Theorem 1 is that we see, in $\mu$, the effect the whole spectrum of the covariance matrix has on the largest eigenvalue of the empirical covariance matrix. This is very different from the classical situation (i.e., $p$ fixed and $n$ goes to $\infty$) where (at least in the real case and when all the eigenvalues of $\Sigma_p$ have multiplicity 1)

$$\sqrt{n}\left(\frac{l_1}{n} - \lambda_1\right) \Rightarrow \mathcal{N}(0, 2\lambda_1^2).$$

(See [2], Theorem 13.5.1.)

In other words, in the classical case, a test based on the largest eigenvalue of the empirical covariance matrix is not sensitive to the whole covariance structure but just to the value of the true largest eigenvalue. Under the asymptotics we are considering, such a test does—implicitly—take into account the whole structure of the spectrum. This is of course very interesting, for instance, for tests of sphericity.

4.2.3. *On Theorem 2 and other random matrices of interest.* The joint distribution of the eigenvalues of other random matrices with complex Gaussian entries is also known. A good reference is, for instance, [23], Section 8. Note that they all involve so-called hypergeometric functions of two matrix arguments. An interesting characteristic of these functions (which since we are dealing with complex entries have to do with the unitary group) is that they have representations in terms of determinants. We refer to, for instance, Section 4 in [20] for explanations and in particular to their Theorem 4.2.

The Harish–Chandra–Itzykson–Zuber formula, which is a preliminary to the proof of Theorem 2, is a subcase of Theorem 4.2 of [20], specialized to the case of the exponential function. A natural question is therefore to know whether one can obtain the same type of representation as the one obtained by Baik–Ben Arous–Johansson–Péché in Theorem 2 in the case of the more general distributions described in [23], Section 8.

In other respects, let us also note the interesting recent developments found in [9] and [15] concerning problems that are close to the one we studied. For more statistical considerations, in the case of spiked models, see [29].

4.3. *Concluding remarks.* The problem of convergence of the joint distribution of the $k$-largest eigenvalues of $X^*X$ requires other tools than the one we discussed in the main body of the paper. We therefore refer the reader to Appendix A.6 for the proof of Corollary 2. In this subsection, we will keep discussing some properties of the largest eigenvalue of $X^*X$.

4.3.1. *Convergence in probability and a.s. convergence.* In this part of the text only, we highlight the fact that $\mu$ depends on $\Sigma_p$, $n$ and $p$ by calling it $\mu(\Sigma_p, n, p)$. Using Slutsky's lemma, it is clear that in the setting of Theorem 1 or 3, for models in $\mathcal{G}$,

$$\frac{l_1}{n} - \mu(\Sigma_p, n, p) \to 0 \qquad \text{in probability.}$$

Since $\mu(\Sigma_p, n, p) > 1/c(\Sigma_p, n, p)$ and $\limsup \lambda_1 c < 1$, we see that $l_1/n$ is always an inconsistent estimator of $\lambda_1$ for models in the class $\mathcal{G}$. Note that Theorem 1 allows us to quantify $(l_1/n) - \lambda_1$ and explore how this quantity is affected by changes in $\Sigma_p$, $n$ and $p$. In particular, elementary computations show that, at $\Sigma_p$ fixed, $(l_1/n) - \lambda_1$ is, unsurprisingly, a decreasing function of $n/p$. We explain in Appendix A.5 that, as announced in Corollary 1, through Theorem 1 or 3 and concentration of measure arguments, we can show that, when the theorems apply,

$$\frac{l_1}{n} - \mu(\Sigma_p, n, p) \to 0 \qquad \text{a.s.}$$

In other respects, we have the following fact.

FACT 2. *Let* $\{Y_{i,j}\}$ *be i.i.d. random variables, real or complex, with* $\mathbf{E}(Y_{i,j}) = 0$, $\mathbf{E}(|Y_{i,j}|^2) = 1$ *and* $\mathbf{E}(|Y_{i,j}|^4) < \infty$. *Let the* $n \times p$ *matrix* $X$ *be such that* $X = Y\Sigma_p^{1/2}$, *where* $Y$ *is an* $n \times p$ *matrix whose entries are the* $Y_{i,j}$. *Suppose the model* $\{(\Sigma_p, n, p)\}$ *is in* $\mathcal{G}$ *and moreover* $H_p \Rightarrow H_\infty$, $n/p \to \rho$ *and* $\lambda_1(\Sigma_p) \to \sup support H_\infty$. *Then*

$$\frac{l_1}{n} - \mu(\Sigma_p, n, p) \to 0 \qquad a.s.,$$

*where* $\mu(\Sigma_p, n, p)$ *is defined in (2).*

It is a simple consequence of Theorem 1.1 and its corollary in [4], once we realize that all the limiting quantities involved in that statement are independent of the distributional assumptions made on the $Y_i$'s. Hence the limit in the case of complex Wishart matrices is the same as the limit in the "general" situation. In particular, this covers the case of real Wishart matrices, that is, data matrices with real normal entries.

4.3.2. *Some simulations.* It was remarked in [25] that the quality of the Tracy–Widom approximation to the marginal distribution of $l_1$ is very good, especially in the right tail of the distribution. This is one of the remarkable properties of this approximation. We refer to [25], Table 1, page 302 for examples. As an aside, we note that the simulation mentioned there was not done with complex Wishart matrices, but rather with real random variables. Nevertheless the same observations hold in the case of complex Wishart matrices with Id covariance. We refer to [11] for theoretical considerations that help understand why this is happening and some simulations in the complex Wishart case.

We made a few simulations to show that the same phenomenon seems to occur in the more involved setting we treat in this paper. Note that numerically solving (1), (2), (3) and getting approximations for $c$, $\mu$ and $\sigma$ takes a fraction of a second on modern computers. We present some results of our experiments in this discussion. See Tables 1 and 2.

We also did some simulations with real Wishart matrices instead of complex ones. In the setting of Theorem 1, we obtained a very reasonable agreement between the empirical distribution of $l_1(X'X)$ and a Tracy–Widom approximation, this time using the Tracy–Widom law appearing in the study of GOE, but keeping the $c$, $\mu$ and $\sigma$ obtained in Theorem 1.

We would finally like to point out that Theorem 1 is essentially explicit if one has access to a computer. Then the eigenvalues of $\Sigma_p$ are numerically computable and so are $c$, $\mu$ and $\sigma$. This is of course a very important property for the relevance of the theorem in applications.

TABLE 1
*Toeplitz covariance matrix example*

| **TW quantiles** | **TW** | **100 × 50** | **400 × 50** | **2* SE** |
|---|---|---|---|---|
| −3.73 | 0.01 | 0.004 | 0.007 | 0.002 |
| −3.20 | 0.05 | 0.033 | 0.041 | 0.004 |
| −2.90 | 0.10 | 0.072 | 0.089 | 0.006 |
| −2.27 | 0.30 | 0.269 | 0.292 | 0.009 |
| −1.81 | 0.50 | 0.479 | 0.497 | 0.010 |
| −1.33 | 0.70 | 0.691 | 0.702 | 0.009 |
| −0.60 | 0.90 | 0.901 | 0.908 | 0.006 |
| −0.23 | 0.95 | 0.953 | 0.956 | 0.004 |
| 0.48 | 0.99 | 0.991 | 0.992 | 0.002 |

The simulation mechanism was as follows. We generated 10,000 random matrices $X$ of size $n \times p$ (using Matlab). The rows of these matrices were i.i.d. $\mathcal{N}_{\mathbb{C}}(0, \Sigma_p)$. For each individual $X$, we computed $l_1(X^*X)/n$ and recentered and rescaled it according to Theorem 1. After simulating 10,000 times we obtained an empirical distribution $\hat{F}$ for $(l_1 - \mu)/(\sigma n^{1/3})$. The columns of the matrix show the value of $\hat{F}$ at the quantiles of the Tracy–Widom distribution (courtesy of Professor Iain Johnstone), given in the leftmost column. If the approximation were "perfect," the third and fourth columns would be equal to the second one.

Here we picked $\Sigma_p = \text{Toeplitz}(1, 0.2, 0.3)$, $p = 50$. For the first column, $n = 100$, $\mu = 3.7297$, $\sigma = 3.9271$. For the second column, $n = 400$, $\mu = 2.6559$, $\sigma = 4.4288$.

TABLE 2
*Sum of atoms example*

| **TW quantiles** | **TW** | **100 × 50** | **400 × 50** | **2* SE** |
|---|---|---|---|---|
| −3.73 | 0.01 | 0.006 | 0.008 | 0.002 |
| −3.20 | 0.05 | 0.036 | 0.045 | 0.004 |
| −2.90 | 0.10 | 0.079 | 0.092 | 0.006 |
| −2.27 | 0.30 | 0.283 | 0.292 | 0.009 |
| −1.81 | 0.50 | 0.490 | 0.496 | 0.010 |
| −1.33 | 0.70 | 0.700 | 0.697 | 0.009 |
| −0.60 | 0.90 | 0.896 | 0.902 | 0.006 |
| −0.23 | 0.95 | 0.949 | 0.951 | 0.004 |
| 0.48 | 0.99 | 0.991 | 0.992 | 0.002 |

The simulation mechanism is similar to the one described previously. We again did 10,000 repetitions of the experiment.

Here $p = 100$. $\Sigma_p$ has $\lambda_1 = \cdots = \lambda_{30} = 10$ and $\lambda_{31} = \cdots = \lambda_{100} = 1$. In the case $n = 100$, $\mu = 24.703$ and $\sigma = 21.871$. In the case $n = 400$, $\mu = 16.417$ and $\sigma = 21.257$.

# APPENDIX A

## A.1. Uniform control of $\Re(f(d(\Sigma_p, n, p)))$ and $\Re(-f(e(\Sigma_p, n, p)))$.

A.1.1. *Case of* $\Re(f(d(\Sigma_p, n, p)))$. Recall that we want to show that $\Re(f(d(\Sigma_p, n, p)))$ is bounded below so as to guarantee that $R_1$, which appears in Lemma 1, is uniformly bounded. In the notation of Section 3.4.1 this is equivalent to showing that

$$m(2(1/\alpha_1 - 1)) \qquad \text{is bounded below.}$$

We clearly have

$$\begin{aligned} &m(2(1/\alpha_1 - 1)) \\ &\quad \geq -\mu c(1/\alpha_1 - q/c) \\ &\qquad - \frac{1}{2\gamma^2}\int \log((1-\alpha)^2 - 2(1/\alpha_1 - 1)\alpha(1-\alpha) \\ &\qquad\qquad + \alpha^2 4(1/\alpha_1 - 1)^2)\, dH_p(\lambda). \end{aligned}$$

It is clear that $((1-\alpha)^2 - 2(1/\alpha_1 - 1)\alpha(1-\alpha) + \alpha^2 4(1/\alpha_1 - 1)^2) \leq 1 + \alpha_1^2 4(1/\alpha_1 - 1)^2$. Note that this quantity is bounded. Note also that the same is true of $1/\gamma^2$, $\mu$ (because $\limsup \alpha_1 < 1$), and hence $-\mu c \alpha_1$ is bounded below. All these arguments together show that $m(2(1/\alpha_1 - 1))$ is bounded below and we have the control we need.

A.1.2. *Case of* $\Re(-f(e(\Sigma_p, n, p)))$. We now want to show that the quantity $\Re(-f(e(\Sigma_p, n, p)))$ is bounded below so that $R_2$ (see Lemma 2) is bounded. In the notation of Section 3.5.1, we need to show $\phi_1(2)$ is bounded below. This quantity is equal to

$$\phi_1(2) = -\frac{1}{2}\log(3) + \frac{1}{2\gamma^2}\int \log(1+\alpha^2)\, dH_p(\lambda) + \frac{1}{2}\left(\frac{1}{\gamma^2} - 1\right)\log(c^2) + \mu q.$$

Now $|\log(c^2)|$ is bounded, since $c$ is bounded away from 0 (see Appendix B.1.1) and $c < 1/\lambda_1 \leq 1/\lambda_p$ and we assume that $\liminf \lambda_p > 0$. Therefore, $\phi_1(2)$ is bounded below and the needed control is shown.

**A.2. About $n^{1/3}$ scaling and its connection to having a saddle point of order 2.** We want to stress that $n^{1/3}$ is the "natural" rate for convergence to Tracy–Widom limits, as there is a connection between Airy functions and saddle points of order 2. The few lines that follow are the natural heuristic explanations of steepest descent analysis. Similar arguments are given after (112) in [6] but we thought it was important to mention them again (and highlight the key parts) since they intuitively explain the connection between having $f''(c) = 0$, an $n^{1/3}$ scaling in Theorem 1 and a Tracy–Widom limit. In another context, the same connections were observed in [17].

Recall that

$$\begin{aligned} A_{n,p}(x) &= -\frac{n\sigma_{n,p}}{2\pi i}\int_\Gamma e^{-n\sigma_{n,p}x(z-q)}e^{-n\mu_{n,p}(z-q)}z^n\prod_{k=1}^p\frac{1}{1-z\lambda_k}\,dz \\ &= -\frac{n\sigma_{n,p}}{2\pi i}\int_\Gamma e^{-n\sigma_{n,p}x(z-q)}e^{nf(z)}\,dz. \end{aligned}$$

Now because $f'(c) = f''(c) = 0$, we have around $c$, $f(z) \simeq f(c) + f^{(3)}(c)(z-c)^3/6$. The point of the steepest descent analysis is to show that we then have (rigorously and up to precision we control)

$$A_{n,p}(x) \simeq -\frac{n\sigma_{n,p}}{2\pi i}\int_\Gamma e^{-n\sigma_{n,p}x(z-q)}e^{n(f(c)+f^{(3)}(c)(z-c)^3/6)}\,dz.$$

Since we picked $f^{(3)}(c) = 2\sigma_{n,p}^3 n^2$, we have

$$\begin{aligned} &e^{-n\sigma_{n,p}x(z-q)}e^{n(f(c)+f^{(3)}(c)(z-c)^3/6)} \\ &\qquad = e^{-n\sigma_{n,p}x(c-q)}e^{nf(c)}\exp\biggl(-xn\sigma_{n,p}(z-c) + \frac{n^3\sigma_{n,p}^3(z-c)^3}{3}\biggr). \end{aligned}$$

A key point is that the Airy function can be written for an appropriately chosen contour $\mathcal{L}$ (see, e.g., [28], page 53):

$$\mathrm{Ai}(x) = \frac{1}{2\pi i}\int_{\mathcal{L}}\exp\biggl(-xv + \frac{v^3}{3}\biggr)\,dv.$$

So the change of variable $a = \tau(z) = n\sigma_{n,p}(z-c)$ becomes natural and our integral can be rewritten as

$$A_{n,p}(x) \simeq -\frac{e^{nf(c)}}{2\pi i}e^{-n\sigma_{n,p}x(c-q)}\int_{\tau(\Gamma)}\exp\biggl(-xa + \frac{a^3}{3}\biggr)\,da.$$

Picking $q = c - \frac{\varepsilon}{n\sigma_{n,p}}$ as in (14), we finally see that

$$A_{n,p}(x) \simeq -\frac{e^{nf(c)}}{2\pi i}e^{-\varepsilon x}\int_{\tau(\Gamma)}\exp\biggl(-xa + \frac{a^3}{3}\biggr)\,da,$$

and the problem is finally to pick a "good" $\Gamma$ on which to analyze $f$ and such that $\tau(\Gamma)$ is an appropriate path from the point of view of the definition of the Airy function. What is on the right-hand side now looks very much like $e^{-\varepsilon x}\mathrm{Ai}(x)/\kappa_{n,p}$ in the notation of Proposition 1. (The minus is of course not a problem since it is an artifact of the orientation of our contours.)

**A.3. Proof of Corollary 3 and examples of models belonging to $\mathcal{G}$.** In this subsection, we show that under assumptions that are both reasonable from an applications standpoint and relatively easy to check, Theorem 1 holds. As in Theorem 1 we assume that $1 \le n/p$, $\limsup n/p < \infty$, $\limsup \lambda_1(\Sigma_p) < \infty$ and $\liminf \lambda_p(\Sigma_p) > 0$. As seen in Appendix B.1.1, these three assumptions imply that $\liminf c > 0$ and $\liminf \lambda_p c > 0$. Our only problem will therefore be to check that

$$\limsup \lambda_1 c < 1.$$

We will use the notation $\alpha = \lambda c$, $\alpha_1 = \lambda_1 c$ and $\gamma^2 = n/p$.

We consider covariance matrices $\Sigma_p$ with spectral distribution $H_p$. We will treat two cases: when $H_p$ has an atom of mass $\nu(p)$ at $\lambda_1$, and the case where $H_p$ weakly converges to a limit and the endpoints of its support converge to the endpoints of the limiting support.

A.3.1. *Case of $H_p$ having an atom of mass $\nu(p)$ at $\lambda_1$.* We assume that $\liminf \nu(p) > 0$. Note that $\lambda_1(\Sigma_p)$ can vary in the analysis that follows. It just needs to be bounded. Since

$$\gamma^2 = \int \frac{\alpha^2}{(1-\alpha)^2} \, dH_p \ge \nu(p) \left( \frac{\alpha_1}{1-\alpha_1} \right)^2,$$

simple algebra shows that

$$\alpha_1 \le \frac{1}{\sqrt{\nu(p)}/\gamma + 1}.$$

Recall that we assume that $\liminf \nu(p) > 0$ and $\limsup n/p < \infty$, so it is clear that $\liminf \sqrt{\nu(p)}/\gamma > 0$ and hence

$$\limsup \alpha_1 < 1$$

in this situation. Therefore Theorem 1 applies.

A.3.2. *Case of weak convergence of $H_p$ with conditions on its support.* We assume that:

1. $H_p \Rightarrow H_\infty$ in the usual weak convergence sense.
2. $\lambda_1(\Sigma_p) \to \sup \text{support} H_\infty \triangleq \lambda_1(\infty)$. We assume that $\limsup \lambda_1(\Sigma_p) < \infty$, so $\lambda_1(\infty) < \infty$.
3. $\lambda_p(\Sigma_p) \to \inf \text{support} H_\infty \triangleq \lambda_\infty(\infty)$ and $\lambda_\infty(\infty) > 0$.
4. In a (left) neighborhood of $\lambda_1(\infty)$, $dH_\infty(\lambda)$ has the property that $dH_\infty(\lambda) \ge B(\lambda_1(\infty) - \lambda)\, d\lambda$, for some $B > 0$.

Hence the only property we have to show is that $\limsup \lambda_1 c < 1$.

Now suppose $H_p \Rightarrow H_\infty$ and $H_\infty$ has a density. Note that for all $x \in [0, -\chi + 1/\lambda_1(\infty)]$, for some $\chi > 0$,

$$t(\lambda) = \left(\frac{\lambda x}{1 - \lambda x}\right)^2$$

is a bounded continuous function of $\lambda$, for $\lambda \in [0, \chi/2 + \lambda_1(\infty)]$. Hence, if we denote

$$f_p(x) = \int \frac{(\lambda x)^2}{(1 - \lambda x)^2} \, dH_p(\lambda),$$

we have

$$f_p(x) \to f_\infty(x) = \int \frac{(\lambda x)^2}{(1 - \lambda x)^2} \, dH_\infty(\lambda),$$

since for $p$ large enough, both $H_p$ and $H_\infty$ are supported in $[0, \chi/2 + \lambda_1(\infty)]$. Now suppose there exist $B > 0$ and $\lambda_B$ such that $dH_\infty(\lambda)/d\lambda \geq B(\lambda_1(\infty) - \lambda)$ in $[\lambda_B, \lambda_1(\infty)]$. Then of course

$$\begin{aligned} f_\infty(x) = \int \frac{(\lambda x)^2}{(1 - \lambda x)^2} \, dH_\infty(\lambda) &\geq \int_{\lambda_B}^{\lambda_1(\infty)} \frac{\lambda^2}{(1/x - \lambda)^2} \, dH_\infty(\lambda) \\ &\geq B \lambda_B^2 \int_{\lambda_B}^{\lambda_1(\infty)} \frac{\lambda_1(\infty) - \lambda}{(1/x - \lambda)^2} \, d\lambda. \end{aligned}$$

Note that

$$\upsilon(x) \triangleq \int_{\lambda_B}^{\lambda_1(\infty)} \frac{\lambda_1(\infty) - \lambda}{(1/x - \lambda)^2} \, d\lambda = \log\left(\frac{1/x - \lambda_B}{1/x - \lambda_1(\infty)}\right) - 1 + \frac{1/x - \lambda_1(\infty)}{1/x - \lambda_B}.$$

Elementary manipulations show that $\upsilon$ is a continuous, increasing function of $x$ on $(0, 1/\lambda_1(\infty))$, going from 0 to $\infty$.

The definition of $f_\infty$ implies that it is a continuous, nondecreasing function of $x$ on the interval $[0, 1/\lambda_1(\infty))$. Since

$$f_\infty(x) \geq B \lambda_B^2 \upsilon(x),$$

we see that $\lim_{x \to \lambda_1(\infty)} f_\infty(x) = +\infty$. Therefore, we can find $b$ such that $f_\infty(b) = 2(1 + \sup n/p)$ and $b$ is bounded away from $1/\lambda_1(\infty)$.

Now recall that $f_p$ is a continuous, increasing function of $x$. Since $f_p(c) = n/p$, when $p$ is large enough, $c \leq b$, since $f_p(b) \to 2(1 + \sup n/p)$. But, for $p$ large enough, $\lambda_1 c \leq \lambda_1 b \to \lambda_1(\infty) b < 1$. Hence $\limsup \lambda_1 c < 1$.

A.3.3. *Some simple examples of matrices for which Theorem 1 applies.* We now justify the claims made after the statement of Corollary 3. We assume that $\limsup \lambda_1(\Sigma_p) < \infty$, $\liminf \lambda_p(\Sigma_p) > 0$, $n \geq p$ and $n/p$ is bounded.

*Sums of atoms.* Suppose $\Sigma_p$ has a largest eigenvalue of multiplicity $k(p)$ and that in the models under consideration $\liminf k(p)/p > 0$. Then we just saw that Theorem 1 applies.

*Equally spaced eigenvalues on an interval.* Suppose the covariance matrices $\Sigma_p$ in our models have eigenvalues that are equally spaced on a fixed interval $[\zeta, \xi]$. Suppose also that $n/p$ is bounded. Then it is clear that the conditions under which we worked in Appendix A.3.2 are satisfied, as long as $\zeta > 0$ and $\xi < \infty$. Hence Theorem 1 applies.

A.3.4. *The case of Toeplitz matrices.* Since we are working with covariance matrices, our matrices $\Sigma_p$ have to be symmetric and positive definite. Let us denote the parameters defining the Toeplitz matrix by $a_0, a_1, \ldots$. Not aiming for the greatest generality, we assume that

$$\sum k|a_k| < \infty.$$

Then the function

$$a(\omega) = a_0 + 2\sum_{k=1}^{\infty} a_k \cos(k\omega)$$

is $\mathcal{C}^1$ on $[0, 2\pi]$. Hence it is bounded and continuous. This function plays an important role in the understanding of the limiting distribution of Hermitian Toeplitz matrices. The results concerning Toeplitz matrices we need are very well known and classical. They can be found in [19], Chapter 5, [18], Chapter 4, and [8], Section 5.5.

Let us denote by $F$ the measure defined on the Borel sets of $\mathbb{R}$ by the following relation: if $E \subset \mathbb{R}$ is a Borel set,

$$F(E) = \frac{1}{2\pi}\mathrm{Leb}\{\omega \in [0, 2\pi] : a(\omega) \in E\},$$

where Leb denotes Lebesgue measure.

As before, we call $H_p$ the spectral measure of $\Sigma_p$, which is now a $p \times p$ Toeplitz matrix. We call $\lambda_1(\infty) = \sup \operatorname{support} F$ and $\lambda_\infty(\infty) = \inf \operatorname{support} F$.

Here is a collection of some interesting and relevant properties of symmetric Toeplitz matrices. Since $a$ is bounded on $[0, 2\pi]$, we have, using Corollary 5.12 in [8], $H_p \Rightarrow F$. $a$ is also piecewise continuous, so $\lim \lambda_p(\Sigma_p) \to \lambda_\infty(\infty)$ and $\lim \lambda_1(\Sigma_p) \to \lambda_1(\infty)$, using, for example, Theorem 5.14 in [8] or Lemma 4.2 in [18]. Finally, it is known ([18], Corollary 4.1 or [8], page 141) that if $F$ does not have any atoms, then its cumulative distribution function $D$ satisfies

$$D(x) = F((-\infty, x]) = \frac{1}{2\pi}\int_{a(\omega)\le x} d\omega.$$

Recall our assumptions: $a$ is bounded away from 0, $\mathcal{C}^1$ and its derivative changes sign only a finite number of times in $[0, 2\pi]$. Also, $F$ is assumed to not have atoms. Then we of course have

$$0 < \inf_{[0,2\pi]} a(\omega) = \lambda_\infty(\infty) \quad \text{and} \quad \sup_{[0,2\pi]} a(\omega) = \lambda_1(\infty) < \infty.$$

Also, we can split $[0, 2\pi]$ into, say, $m$ intervals where $a$ is monotonic. Calling the intervals $I_k$ and their endpoints $p_k$ (with $I_1 < I_2 < \cdots$ and $I_1 = [p_1, p_2]$), we have

$$D(x) = \sum_{k=1}^{m} \frac{1}{2\pi} \int_{\omega \in I_k \,:\, a(\omega) \leq x} d\omega.$$

The function $a$ is invertible on $I_k$. Also, $\lambda_1(\infty)$ is reached and so there is at least one $k$, say $k_0$, for which $a(p_{k_0+1}) = \lambda_1(\infty)$. Further, we can assume without loss of generality that $a$ is nondecreasing on $I_{k_0}$. We call $a_{k_0}$ the restriction of $a$ to $I_{k_0}$. $a_{k_0}$ is an invertible function. Now, assuming that $a(p_{k_0+1}) \geq x \geq a(p_{k_0})$, we have

$$D_{k_0}(x) = \int_{\omega \in I_{k_0} \,:\, a(\omega) \leq x} d\omega = a_{k_0}^{-1}(x) - p_{k_0}.$$

Since $a_{k_0}$ is $\mathcal{C}^1$, $D_{k_0}$ has a derivative in $(a(p_{k_0}), a(p_{k_0+1}))$ and we have

$$D'_{k_0}(x) = \frac{1}{a'_{k_0}(a_{k_0}^{-1}(x))}.$$

We immediately see that on this interval

$$D'_{k_0}(x) \geq \frac{1}{\sup_{[0,2\pi]} |a'(\omega)|} > 0$$

since $a$ is $\mathcal{C}^1$.

Hence, after we rewrite $D$ as a sum of $D_k$'s, we see that under our assumptions $D$ has a density except at a finite number of points where the derivative of $a$ changes sign. The density tends to $\infty$ at these points. So the assumptions put forth in Appendix A.3.2 hold and Theorem 1 applies to the class of Toeplitz covariance matrices we considered.

In general, if $a$ is a Lebesgue integrable function on $(-\pi, \pi)$ whose Fourier coefficient coincides with the $a_i$'s, and if $\operatorname{ess\,sup} a = M_a < \infty$ and $\operatorname{ess\,inf} a = m_a > 0$, Theorem 1 holds for such a Toeplitz matrix if

$$T(x) = \int_{-\pi}^{\pi} \left( \frac{a(u)x}{1 - a(u)x} \right)^2 du$$

is a continuous function of $x$ on $[0, 1/M_a)$ that is increasing and tends to $\infty$ as $x \to 1/M_a$. (Note that since $a \geq 0$ a.e., $T$ is nondecreasing in $x$.) This is a simple consequence of the so-called First (or Weak) Szegö limit theorem (see [19], pages 64–65) and of the fact that the eigenvalues of the corresponding (truncated) Toeplitz matrices are between $m_a$ and $M_a$ in this situation.

**A.4. Justification of results for spiked models with a small spike.** Here we are considering "spiked" models of covariance. Namely, we start with a model $\{\Sigma_p, n, p\}_{n,p\in\mathbb{N}}$ that is in $\mathcal{G}$. In other words, Theorem 1 applies to this model. When we say that we are considering the spiked version of this model, we mean that we are now focusing on data matrices $X$ that are $n \times (p+k)$, and $X_i \overset{\text{i.i.d.}}{\sim} \mathcal{N}_{\mathbb{C}}(0, \widetilde{\Sigma}_{p+k})$, where $k(p) < K$, $K \in \mathbb{N}$, and we add to $\{\lambda_1(\Sigma_p), \dots, \lambda_p(\Sigma_p)\}$ $k$ eigenvalues larger than $\lambda_1(\Sigma_p)$. In other words, $\lambda_{k+1}(\widetilde{\Sigma}_{p+k}) = \lambda_1(\Sigma_p)$.

A.4.1. *Proof of Fact 1.* The statement we want to prove is the following: *In the "spiked" situation described above, if there exists $\chi > 0$ such that*

$$\lambda_1(\widetilde{\Sigma}_{p+k}) < -\chi + 1/c(\Sigma_p, n, p),$$

*Theorem 1 applies to $\widetilde{\Sigma}_{p+k}$.*

Proof. In order to simplify the notation we will use in this proof the shortcuts

$$\begin{aligned} \tilde{c} &\triangleq c(\widetilde{\Sigma}_{p+k}, n, p+k), \\ c &\triangleq c(\Sigma_p, n, p), \\ \tilde{\lambda}_1 &\triangleq \lambda_1(\widetilde{\Sigma}_{p+k}). \end{aligned}$$

It is clear that the only thing we have to check is that

$$\limsup \lambda_1(\widetilde{\Sigma}_{p+k}) c(\widetilde{\Sigma}_{p+k}, n, p+k) < 1.$$

We of course have $c < 1/\tilde{\lambda}_1$. Now let us call

$$\rho(x) = \int \left(\frac{\lambda x}{1-\lambda x}\right)^2 d\widetilde{H}_{p+k} \qquad \text{where } \rho \text{ is defined on } [0, 1/\tilde{\lambda}_1).$$

The equation that defines $\tilde{c}$ is

$$\rho(\tilde{c}) = \frac{n}{p+k} \qquad \text{with } \tilde{c} \in [0, 1/\tilde{\lambda}_1).$$

We have seen that $\rho$ is an increasing function of $x$. Now since $c(\Sigma_p, n, p) < 1/\tilde{\lambda}_1$, we can compute $\rho(c)$. Note that we have, if we denote by $\tilde{\lambda}_i$'s the eigenvalues we have added to $\Sigma_p$ to create $\widetilde{\Sigma}_{p+k}$,

$$\rho(x) = \frac{1}{p+k} \sum_{j=1}^{k} \left(\frac{\tilde{\lambda}_j x}{1-\tilde{\lambda}_j x}\right)^2 + \frac{p}{p+k} \int \left(\frac{\lambda x}{1-\lambda x}\right)^2 dH_p(\lambda).$$

Now recall that by definition,

$$\int \left(\frac{\lambda c}{1-\lambda c}\right)^2 dH_p(\lambda) = \frac{n}{p}.$$

Hence

$$\rho(c) = \frac{1}{p+k} \sum_{j=1}^{k} \left( \frac{\tilde{\lambda}_j c}{1 - \tilde{\lambda}_j c} \right)^2 + \frac{n}{p+k} > \frac{n}{p+k}.$$

Since $\rho$ is an increasing function of $x$ this implies that

$$\tilde{c} < c$$

and therefore

$$\tilde{\lambda}_1 \tilde{c} < \tilde{\lambda}_1 c < 1 - \chi c.$$

Since $\liminf c > 0$ because $\{\Sigma_p, n, p\} \in \mathcal{G}$, we have shown

$$\limsup \tilde{\lambda}_1 \tilde{c} < 1$$

and Theorem 1 applies to $\widetilde{\Sigma}_{p+k}$. $\square$

**A.5. Issues of convergence in probability and a.s. convergence.** We will explain in this subsection why, when Theorem 1 or 3 applies, we have

$$\frac{l_1}{n} - \mu(\Sigma_p, n, p) \to 0 \qquad \text{in probability}$$

and

$$\frac{l_1}{n} - \mu(\Sigma_p, n, p) \to 0 \qquad \text{a.s.}$$

The convergence in probability part is an immediate application of Slutsky's lemma (see [39], Lemma 2.8), so we will not belabor this point. The only thing we have to show is therefore the almost sure convergence part. We use concentration of measure arguments to show that $l_1/n - \mu(\Sigma_p, n, p) \to 0$ a.s.

FACT. *If Theorem 1 or Theorem 3 applies,*

$$\frac{l_1}{n} - \mu(\Sigma_p, n, p) \to 0 \qquad \textit{a.s.}$$

PROOF. Let us first recall that the application that takes a matrix $M$ and returns its ordered singular values is 1-Lipschitz with respect to Euclidean norms (see, e.g., statement 7.3.8 in [22]). In other words, if we call $\{\sigma_i\}$ and $\{\tau_i\}$ the ordered singular values of two $n \times p$ matrices $A$ and $B$, we have

$$\sum_{k=1}^{p} (\sigma_k - \tau_k)^2 \le \sum_{i,j} |a_{i,j} - b_{i,j}|^2.$$

In particular, that shows that the application that takes a vector of dimension $2np$, turns it into matrices $M$ and $N$ and returns the ordered singular values of $M+iN$ is 1-Lipschitz with respect to Euclidean norms. So is any 1-Lipschitz (for Euclidean norms) $\mathbb{R}^p \to \mathbb{R}$ function of the ordered singular values, and in particular the projection that returns $\sigma_1$ from $(\sigma_1, \sigma_2, \dots, \sigma_p)$. Hence, by a fairly standard concentration of measure argument (see, e.g., [13], pages 34–38, and references therein, especially [21]), if we call $m(\Sigma_p, n, p)$ a median of $s_1 = \sqrt{l_1(X^*X)/n}$, and $\bar{\lambda}_1 = \sup \lambda_1(\Sigma_p) < \infty$, we have, in the setting of Theorem 1 or 3,

$$\forall r > 0, \qquad P(|s_1 - m(\Sigma_p, n, p)| > r) \leq 2\exp(-nr^2/\bar{\lambda}_1).$$

Note that the fact that the rows of our matrices are $\mathcal{N}_{\mathbb{C}}(0, \Sigma_p)$ plays a crucial role here, for we know the concentration function of Gaussian random variables and we also know that it has the so-called dimension-free concentration property. We refer the reader to [26], page 99, for more information about it. Let us just say that, for quite general distributions, the interplay between log-Sobolev inequalities and concentration of product measures is the gist of the argument that leads to the previous inequality. In the Gaussian case, we can also use the fact that the joint distribution of the entries of the $2np$ vector has a density of the type $\exp(-U)$ with $\text{Hessian}(U) \geq 2\text{Id}/\lambda_1$. [Recall that since we are working with complex standard entries, the rows of $M$ and $N$ are i.i.d. $\mathcal{N}(0, \Sigma_p/2)$.] Hence Theorem 2.8 in [26] applies and the concentration function for this measure is $\exp(-r^2/\lambda_1)$.

Combining it with the first Borel–Cantelli lemma (recall that $n$ is going to $\infty$), we see that

$$s_1 - m(\Sigma_p, n, p) \to 0 \qquad \text{a.s.}$$

Since we know that $\mu(\Sigma_p, n, p)$ is uniformly bounded when Theorem 1 or Theorem 3 applies, we conclude that in this situation $\sqrt{\mu(\Sigma_p, n, p)}$ is, too.

Now because $s_1 - m(\Sigma_p, n, p) \to 0$ a.s., $s_1 - m(\Sigma_p, n, p) \to 0$ in probability. But we also know that $l_1/n - \mu(\Sigma_p, n, p) \to 0$ in probability. Therefore, $\sqrt{l_1/n} - \sqrt{\mu(\Sigma_p, n, p)} \to 0$ in probability, because, for instance, $\mu(\Sigma_p, n, p)$ is bounded below. And so $m(\Sigma_p, n, p) - \sqrt{\mu(\Sigma_p, n, p)} \to 0$.

Hence, there exists $K > 0$ such that $0 \leq s_1 \leq K$ a.s. Hence $s_1$ is (a.s.) uniformly bounded. So is $m(\Sigma_p, n, p)$ and hence

$$\begin{aligned}(s_1 - m(\Sigma_p, n, p))(s_1 + m(\Sigma_p, n, p)) &= s_1^2 - m(\Sigma_p, n, p)^2 \\ &= \frac{l_1}{n} - m(\Sigma_p, n, p)^2 \to 0 \qquad \text{a.s.}\end{aligned}$$

We know that $m(\Sigma_p, n, p)^2 - \mu(\Sigma_p, n, p) \to 0$, because $\mu(\Sigma_p, n, p)$ is bounded above, so we have shown

$$\frac{l_1}{n} - \mu(\Sigma_p, n, p) \to 0 \qquad \text{a.s.} \qquad \square$$

**A.6. Determinantal character of the point process and consequences.** In this section, we explain that when viewed as a point process on the real line, the eigenvalues of $X^*X$ form a determinantal point process. The main consequence is that when the rows $X_i \overset{\text{i.i.d.}}{\sim} \mathcal{N}(0,\Sigma_p)$ and the covariance models are in the class $\mathcal{G}$, the joint distribution of the $k$-largest eigenvalues of $X^*X$ ($k$ fixed) converges to its Tracy–Widom counterpart.

The fact that the point process is determinantal is an easy consequence of a well-known result that seems to first appear in [7], Section 2, and that directly applies to the situation we are considering, given the form of the density of the eigenvalues of $X^*X$, when $X_i$ are $n$ i.i.d. $\mathcal{N}_{\mathbb{C}}(0,\Sigma_p)$. Proposition 2.1 in [6] shows that the kernel of this determinantal point process is $K_{n,p}$, where $K_{n,p}$ is defined in (4). We can now turn to the issue of the convergence of the joint distribution.

A.6.1. *Convergence of the joint distribution.* Let $B_j$ be disjoint, bounded below Borel sets of $\mathbb{R}$ and let $N_{B_j}$ denote the number of eigenvalues of $X^*X$ that are in $B_j$. As explained in Theorem 2 in [33] (see also (2.44) in [34]), the generating function of the probability distribution of $N_{B_j}$ can be written as the determinant of an operator. In our case, if we call $L=\sum_{j=1}^k (z_j-1)1_{B_j}$, we have

$$\mathbf{E}\left(\prod_{j=1}^k z_j^{N_{B_j}}\right)=\det(\mathrm{Id}+K_{n,p}L).$$

Using Lemma 2 in [34], if we can show that $\det(\mathrm{Id}+A_{n,p}B_{n,p}L)\to\det(\mathrm{Id}+\mathrm{Ai}^2L)$, we will have shown the convergence of the joint distribution of the $k$-largest eigenvalues of $X^*X$ (properly recentered and rescaled) to their Tracy–Widom counterpart. (The argument is similar to the one given in the proof of Theorem 1, pages 1047–1048 in [34].)

Now recall that we showed that $A_{n,p}B_{n,p}\to E\mathrm{Ai}^2E^{-1}$ in trace class norm in the notation of Lemma 3. Our only problem is therefore to show that $\det(\mathrm{Id}+E\mathrm{Ai}^2E^{-1}L)=\det(\mathrm{Id}+\mathrm{Ai}^2L)$. Note that since $L$ and $E^{-1}$ are multiplication operators, they commute. Also, recall that $E\mathrm{Ai}$ and $\mathrm{Ai}E^{-1}$ are Hilbert–Schmidt operators. Since $L$ is bounded, $\mathrm{Ai}LE^{-1}$ is Hilbert–Schmidt. Recall also that for Hilbert–Schmidt operators $F$ and $G$, $\det(\mathrm{Id}+FG)=\det(\mathrm{Id}+GF)$. Hence,

$$\begin{aligned}\det(\mathrm{Id}+E\mathrm{Ai}^2E^{-1}L)&=\det(\mathrm{Id}+E\mathrm{Ai}^2LE^{-1})=\det(\mathrm{Id}+(\mathrm{Ai}LE^{-1})(E\mathrm{Ai}))\\&=\det(\mathrm{Id}+\mathrm{Ai}L\mathrm{Ai})=\det(\mathrm{Id}+\mathrm{Ai}^2L).\end{aligned}$$

We refer to [35] and [10] for information about the limiting distributions of $l_2,\dots,l_k$.

## APPENDIX B: CONVERGENCE OF OPERATORS

In this section, we will prove Proposition 1 (which deals with the convergence of $A_{n,p}$) and sketch the proof of Proposition 2 (which does the same thing for $B_{n,p}$). The method of proof is similar to what is done for the proof of Proposition 3.1 in [6]. It might look a little simpler because we worked in the beginning of this paper with more complicated functions $f$ than [6] did. So, from the point of view of this analysis, the efforts are in some sense balanced differently.

At this point, what we have to do is adapt the proof found in [6]; the difficult conceptual and technical problems we had to solve that required fresh ideas and a new look are found earlier in the paper. Now we principally need to rephrase parts of the work of [6] in a more general context, once the gist of the argument is understood in this general context. Note that our paths are slightly different from theirs, and there are a few other things to check. (In particular, we state Proposition 1 with an $\exp(-\varepsilon s/2)$ in the upper bound, independent of the interval $[-s_0, \infty)$ on which we work. We need to show that one can adapt the proof given in [6] to do this and not have $b(s_0)$, possibly dependent on $s_0$, instead of $\varepsilon$.)

We decided to include the full proof for three reasons. A sequence of references to various equations in [6] and modifications to make to those would have made for a very difficult reading. It would also have assumed that the reader had an enormous familiarity with [6]. So we decided to include this analysis for the convenience of the reader. Also, given the somewhat technical nature of the problem, having a completely spelled out proof reduces considerably the risk of errors.

Nevertheless, because of the length of the proofs, we will only give a complete proof for the convergence of $\kappa_{n,p} A_{n,p}$ to its limit. We will just sketch the corresponding proof for $B_{n,p}/\kappa_{n,p}$.

**B.1. Preliminary remarks.** We first recall the assumptions satisfied by models in $\mathcal{G}$. We assume:

1. $n/p$ is uniformly bounded and greater than or equal to 1.
2. $\limsup \lambda_1(\Sigma_p) < \infty$.
3. $\liminf \lambda_p(\Sigma_p) > 0$.
4. $\limsup \lambda_1 c < 1$.

Recall also that $f$, whose dependence on $(\Sigma_p, n, p)$ we choose to not highlight, is defined as

$$f(z) = -\mu(z - q) + \log(z) - \frac{p}{n} \int \log(1 - z\lambda)\, dH_p(\lambda).$$

Before we explain why the proof of Proposition 3.1 in [6] can be adapted to our problem under these assumptions, we need to show an intermediary result: the fact that under the above assumptions, $\liminf c > 0$.

B.1.1. *About* $\liminf c$. The fact that $\overline{\lambda}_1 \triangleq \limsup \lambda_1 < \infty$ implies that $\liminf c > 0$. As a matter of fact we have

$$1 \leq \frac{n}{p} = \int \left(\frac{\lambda c}{1-\lambda c}\right)^2 dH_p(\lambda) \leq \left(\frac{\lambda_1 c}{1-\lambda_1 c}\right)^2$$

and hence

$$\frac{1}{2\lambda_1} \leq c.$$

This of course implies that

$$\liminf c(\Sigma_p, n, p) > 0$$

since we assume $\limsup \lambda_1 < \infty$.

B.1.2. *Key properties needed for the proof to go through.* As explained in Section 3.6.2, there are four crucial points that will allow us to carry out the proof.

The first one is the fact that the lengths of $\Gamma$ and $\Xi$ are uniformly bounded when the (nonspiked) covariance models are in $\mathcal{G}$. It is clear that this is implied by the condition $\liminf \lambda_p c > 0$ (which is equivalent to $\liminf \lambda_p > 0$, since $c$ is bounded below) and the fact that $\Re(f(d))$ and $\Re(-f(e))$ are bounded below under our assumptions (which implies that $R_1$ and $R_2$ are bounded).

The second very important point is that one needs to be able to find $\delta > 0$ such that

$$\exists \delta > 0, \forall s \qquad |s-c| < \delta \quad \Longrightarrow \quad \frac{|f^{(4)}(s)|}{4!}\delta < \frac{\sigma^3}{6}.$$

The importance of this property will become clear in the proof. Of course, this has to be uniform with respect to our models. In our context, calling

$$\limsup \lambda_1 c = \overline{\alpha}_1 \quad \text{and} \quad \delta = \eta c,$$

it is easy to see that this is implied by

$$\frac{\eta}{4c^3}\left[\frac{1}{(1-\eta)^4} + \frac{p}{n}\left(\frac{\overline{\alpha}_1}{1-(1+\eta)\overline{\alpha}_1}\right)^4\right] < \frac{1}{c^3}$$

or

$$\eta\left[\frac{1}{(1-\eta)^4} + \frac{p}{n}\left(\frac{\overline{\alpha}_1}{1-(1+\eta)\overline{\alpha}_1}\right)^4\right] < 4.$$

Since by assumption $\overline{\alpha}_1 < 1$ and $p/n \leq 1$ it is clear that we can find $\eta > 0$ such that the inequality appearing in the previous display is verified.

Therefore, $\delta = \liminf \eta c$ is bounded away from 0, since $\eta$ and $c$ both are. The assumptions $\limsup \lambda_1 c < 1$, $p/n \leq 1$ and the fact that $\liminf c > 0$ imply that both $\mu$ and $\sigma$ [defined by (2) and (3)] are bounded, which insures that for the same $\delta$

$$\forall s \qquad |s-c|<\delta \quad \Longrightarrow \quad \limsup \sup \frac{|f^{(4)}(s)|}{4!} = \Delta < \infty.$$

Finally, note that since $\overline{\alpha}_1 < 1$, we can guarantee that the $\delta$ we pick is small enough that the disc of center $c$ and radius $\delta$ never encloses $d(\Sigma_p, n, p)$. [For obvious symmetry reasons, it also means that $\delta$ can be chosen small enough that $e(\Sigma_p, n, p)$ is not enclosed either.]

B.1.3. *On Theorem 3.* In this situation, we consider the case where a covariance model $\{\Sigma_p, n, p\}$ in $\mathcal{G}$ is spiked with $k$ eigenvalues at $1/c(\Sigma_p, n, p)$ ($k = 0$ is a possibility).

Using notation similar to what we used earlier, we will need to analyze the function

$$\begin{aligned} A_{\widetilde{\Sigma}_{p+k},n,p+k}(x) &= -\frac{n\sigma_{n,p}}{2\pi i}\int_\Gamma e^{-n\sigma_{n,p}x(z-q)} e^{-n\mu_{n,p}(z-q)} z^n \prod_{k=1}^p \frac{1}{1-z\lambda_k} \frac{c^k}{(c-z)^k}\,dz \\ &= -\frac{n\sigma_{n,p}}{2\pi i}\int_\Gamma e^{-n\sigma_{n,p}x(z-q)} e^{-nf_{n,p}(z)} \frac{c^k}{(c-z)^k}\,dz, \end{aligned}$$

where $f_{n,p}$ is the function that appears in the analysis of $A_{\Sigma_p,n,p}$. We used the index $(n, p)$ to remove any ambiguity.

Similarly, we will have to study

$$\begin{aligned} &B_{\widetilde{\Sigma}_{p+k},n,p+k}(x) \\ &\quad = \frac{n\sigma_{n,p}}{2\pi i}\int_\Xi e^{n\sigma_{n,p}x(z-q)} e^{n\mu_{n,p}(z-q)} \frac{1}{z^n} \prod_{k=1}^p (1-\lambda_k z)(c-z)^k/c^k\,dz. \end{aligned}$$

What we will show is that after proper scaling, these functions converge to limiting functions $\tilde{H}_{\infty,k}$ and $\tilde{J}_{\infty,k}$, defined thereafter (and appearing first in (120) and (122) in [6]).

Theorem 3 is a generalization of Theorem 1.1(a) in [6] in the sense that it shows that the same limiting distributions $F_k$'s appear if we spike the covariance matrix at the "critical" eigenvalue. Note nevertheless that we do not recover exactly the same critical eigenvalue. We could have if we looked at a model of the type $\{\Sigma_{p+k+r}, n, p+k+r\}$, with $r$ eigenvalues such that for some $\chi > 0$, $l_{k+1}, \dots, l_r \in [\chi, 1/c(\Sigma_p, n, p) - \chi]$. This would have added a little bit of technical difficulty to the proof we give later without the benefit of understanding since we already saw that the model $\{\Sigma_{p+r}, p+r, n\}$ (corresponding to $\Sigma$ and those $r$ "extra" eigenvalues) is in $\mathcal{G}$.

B.1.4. *Limiting functions and limiting distributions.* It seemed to us that slight (essentially "cosmetic") modifications of the functions introduced in [6], (120)–(122), were the most natural way to define them, especially when considering existing literature on Tracy–Widom limits. So we call

$$H_{\infty,k}(x) = -\frac{1}{2\pi i}\int_{\Gamma_\infty} \exp(-ax + a^3/3)\frac{da}{a^k}. \tag{15}$$

Here if we call $\varepsilon$ the positive real introduced earlier in the text, $\Gamma_\infty$ goes from $\infty e^{i\pi/3}$ to $\infty e^{-i\pi/3}$, goes through the real axis on the left of 0, stays in the region $\Re(z+\varepsilon) > 0$ and is symmetric about the real axis. It is oriented counterclockwise. In subsequent analysis, we will take $\Gamma_\infty$ to be the union of the straight line $te^{i\pi/3}$, $\infty > t \geq \varepsilon/2$, the arc of circle of center 0 and radius $\varepsilon/2$, for angles $\theta \in [\pi/3, 5\pi/3]$, and the straight line $te^{-i\pi/3}$ for $\varepsilon \leq t < \infty$. Note that when $k = 0$, $H_{\infty,0}(x) = \mathrm{Ai}(x)$.

Similarly, let

$$J_{\infty,k}(x) = \frac{1}{2\pi i}\int_{\Xi_\infty} \exp(ax - a^3/3)a^k\, da. \tag{16}$$

Here, the contour $\Xi_\infty$ is restricted to the region $\Re(z+\varepsilon) < 0$, goes from $\infty e^{-i2\pi/3}$ to $\infty e^{i2\pi/3}$ and is symmetric about the real axis. It is also oriented counterclockwise. In subsequent analysis, we will take it to be the union of the line $te^{-i2\pi/3}$, $3\varepsilon < t < \infty$, the arc of circle, $3\varepsilon e^{i\theta}$, $\theta \in [2\pi/3, 4\pi/3]$ and the line $te^{i2\pi/3}$, $3\varepsilon < t < \infty$.

Note that $\Xi_\infty$ is strictly to the left of $\Gamma_\infty$.

Finally, using (206) in [6], it is clear that $|e^{-\varepsilon x}H_{\infty,k}(x)| \leq Ke^{-\varepsilon x/2}$, for some $K > 0$. Using (205) there, we get similarly that $e^{\varepsilon x}J_{\infty,k}(x) = \mathrm{O}(e^{-\varepsilon x/2})$ on $[s_0,\infty)$, for all $s_0 > \infty$.

Hence, if we call $\tilde{H}_{\infty,k}(x) = e^{-\varepsilon x}H_{\infty,k}(x)$ and $\tilde{J}_{\infty,k}(x) = e^{\varepsilon x}J_{\infty,k}(x)$, we see that the operators on $L^2([0,\infty))$ with kernel $K(x,y) = \tilde{H}_{\infty,k}(x+y+s)$ and $k(x,y) = \tilde{J}_{\infty,k}(x+y+s)$ are Hilbert–Schmidt, for any fixed $s$. Hence their product is trace class.

The cumulative distribution functions $F_k$'s mentioned in Theorem 3 are connected to $\tilde{H}_{\infty,k}$ and $\tilde{J}_{\infty,k}$ in the following manner. If we call, by a slight abuse of notation, $\tilde{H}_{\infty,k}$ the operator with kernel $\tilde{H}_{\infty,k}(x+y+s)$ on $L^2[0,\infty)$ and similarly $\tilde{J}_{\infty,k}$ the operator with kernel $\tilde{J}_{\infty,k}(x+y+s)$ on the same space, then

$$F_k(s) = \det(I - \tilde{H}_{\infty,k}\tilde{J}_{\infty,k}).$$

Note that as explained in [6], this quantity is well defined and independent of $\varepsilon$. $\varepsilon$ is just here to ensure that $\tilde{H}_{\infty,k}$ is Hilbert–Schmidt.

**B.2. Convergence of $A_{n,p}$.** We work in the general case where there is a root of multiplicity $k$ at $c(\Sigma_p, n, p)$. We nevertheless will not highlight this dependence on $k$ to simplify the notation.

Denoting by $f$ the function defined in (11) and corresponding to $\{\Sigma_p, n, p\}$, we call $\kappa_{n,p} = \exp(-nf(c))/(-\sigma n^{1/3} c)^k$ and, since we are in the case where $\sigma_{n,p} = \sigma/n^{-2/3}$,

$$A_{\widetilde{\Sigma}_{p+k},n,p+k}(s) = A_{n,p}(s) = -\frac{n^{1/3}\sigma}{2\pi i}\int_\Gamma e^{-n^{1/3}s(z-q)} e^{nf(z)} \frac{c^k\,dz}{(c-z)^k}.$$

Hence,

$$\mathcal{A}_{n,p} \triangleq \kappa_{n,p} A_{n,p}(s) = -\frac{1}{2\pi i}\frac{1}{(\sigma n^{1/3})^{k-1}}\int_\Gamma e^{-n^{1/3}\sigma s(z-q)} e^{n(f(z)-f(c))}\frac{dz}{(z-c)^k}.$$

The aim of this subsection is to show the following lemma.

LEMMA B.1. *Let $f$ satisfy Lemma 1. Suppose*

$$\exists \delta > 0 : \forall s \qquad |s-c| < \delta \Rightarrow \left|\frac{f^{(4)}(s)}{4!}\right|\delta \le \frac{\sigma^3}{6}.$$

*Then*

$$\forall s_0 \in \mathbb{R}, \exists C(s_0), \exists N_0(s_0) : s > s_0, n > N_0$$
$$\Rightarrow |\kappa_{n,p} A_{n,p}(s) - \exp(-\varepsilon s) H_{\infty,k}(s)| \le \frac{C(s_0)\exp(-\varepsilon s/2)}{n^{1/3}}.$$

We now turn to proving Lemma B.1.

B.2.1. *Notation.* We call $\mathcal{C}$ the circle of center $c$ and radius $\delta$. We call $\mathcal{D}$ the corresponding disc. We split $\Gamma$ into $\Gamma = G^{(1)} \cup G^{(2)}$, where $G^{(1)}$ is the part of $\Gamma$ that is inside $\mathcal{D}$ (see Figure 4). Note that under our assumptions about $\delta$ and $d$, the intersection of $\Gamma$ and $\mathcal{C}$ is on $\Gamma_1 \cup \overline{\Gamma_1}$, that is, on a section of $\Gamma$ where this contour is parametrized as $c + te^{\pm i\pi/3}$, $t \in \mathbb{R}$.

We call $\Gamma_\infty^{(1)}$ the image of $G^{(1)}$ under the map $z \mapsto \sigma n^{1/3}(z-c)$. Of course, everything has been done so that this is a subset of $\Gamma_\infty$. Let us denote $\Gamma_\infty^{(2)} = \Gamma_\infty \setminus \Gamma_\infty^{(1)}$.

Recall that we called $\mathcal{A}_{n,p}(x) = \kappa_{n,p} A_{n,p}(x)$. Let $\mathcal{A}_{n,p}(x) = \mathcal{A}_{n,p}^{(1)}(x) + \mathcal{A}_{n,p}^{(2)}(x)$, where the superscript indicates that $\mathcal{A}_{n,p}^{(i)}(x)$ is the contribution of the integral defining $\mathcal{A}_{n,p}(x)$ over $G^{(i)}$.

We similarly split $H_{\infty,k}$ into $H_{\infty,k} = H_{\infty,k}^{(1)} + H_{\infty,k}^{(2)}$ where now the superscripts refer to the contribution of the integrals over $\Gamma_\infty^{(1)}$ and $\Gamma_\infty^{(2)}$.

B.2.2. *Preliminary computations.* Recall that $\sigma$ is uniformly bounded (away from 0 and $\infty$) for models in the class $\mathcal{G}$. Since we supposed that $|f^{(4)}(s)/4!|\delta \le \sigma^3/6$ and $\delta$ is bounded away from 0, it is clear that there exists $0 < \Delta < \infty$ such that $\sup_{|s-c|\le\delta}|f^{(4)}(s)/4!| \le \Delta$, uniformly for our models.

We now turn to bounding a quantity that is key in the analysis. We have, for any complex number $z$, $|\Re(z)| \le |z|$. Therefore, since $f$ has two 0 derivatives at $c$, we have by Taylor's theorem, for $z$'s such that $|z-c| \le \delta$,

$$\begin{aligned}\left|\Re\left(f(z)-f(c)-\frac{f^{(3)}(c)}{6}(z-c)^3\right)\right| &\le \left|f(z)-f(c)-\frac{f^{(3)}(c)}{6}(z-c)^3\right|,\\ &\le \left(\sup_{|s-c|\le\delta}\frac{|f^{(4)}(s)|}{4!}\right)|z-c|^4,\\ &\le \left(\sup_{|s-c|\le\delta}\frac{|f^{(4)}(s)|}{4!}\right)\delta|z-c|^3,\\ &\le \frac{\sigma^3}{6}|z-c|^3,\end{aligned}$$

because $|z-c| \le \delta$. Recall that $f^{(3)}(c) = 2\sigma^3$. Hence when $z \in \mathcal{D}$,

$$\Re(f(z)-f(c)) \le \frac{\sigma^3}{6}(2\Re((z-c)^3) + |z-c|^3).$$

In particular, when $z \in \mathcal{D}$ and $z = c + te^{\pm i\pi/3}$, $\Re(f(z)) \le f(c) - t^3\sigma^3/6$. Now recall that because $\Re(f)$ is decreasing on $\Gamma_1$ and since $d \notin \mathcal{D}$, $f(d) \le f(c+\delta e^{i\pi/3})$. If $z$ is in $G^{(2)}$, either it is on $\Gamma_1$ or $\Re(z) \ge \Re(d)$. In the latter case, $\Re(f(d)) \ge \Re(f(z))$ because $f$ satisfies Lemma 1. In the former, we can use the fact that $\Re(f(z))$ is decreasing on $\Gamma_1$ to finally get that for $z \in G^{(2)}$, $\Re(f(z)) \le \Re(f(c+\delta e^{i\pi/3})) \le f(c) - \sigma^3\delta^3/6$.

We will now split the analysis into three parts corresponding to different regions of $\Gamma$.

B.2.3. *Behavior of our functions on $G^{(2)}$ and $\Gamma_\infty^{(2)}$.* Let us first focus on $\mathcal{A}_{n,p}^{(2)}(x)$. By definition,

$$\mathcal{A}_{n,p}^{(2)}(x) = -\frac{1}{2\pi i}\frac{1}{(\sigma n^{1/3})^{k-1}}\int_{G^{(2)}} e^{-n^{1/3}\sigma x(z-q)} e^{n(f(z)-f(c))}\frac{dz}{(z-c)^k}.$$

Hence,

$$|\mathcal{A}_{n,p}^{(2)}(x)| \le \frac{1}{2\pi(\sigma n^{1/3})^{k-1}}\int_{G^{(2)}} e^{-n^{1/3}\sigma x\Re(z-q)} e^{n\Re(f(z)-f(c))}\frac{|dz|}{|c-z|^k}.$$

Now on $G^{(2)}$, $|c-z| \ge \delta$ and $\Re(f(z)-f(c)) \le -\sigma^3\delta^3/6$. So we only have to pay close attention to $n^{1/3}\sigma x\Re(z-q)$.

Suppose $x \in [-s_0, \infty)$, with $s_0 > 0$. If $x > 0$, then $-\sigma x \Re(z-q) \leq -\sigma x\delta/2 \leq \sigma s_0 R_1 - \sigma x \delta/2$. If $x < 0$, then $-\sigma x \Re(z-q) \leq \sigma s_0 R_1$, because $\Re(q) > 0$ and we saw that $R_1$ can be chosen to be independent of our models (i.e., uniform with respect to them). We of course also have $-\sigma x \Re(z-q) \leq \sigma s_0 R_1 - \sigma x\delta/2$ when $x < 0$.

So, since the length of $G^{(2)}$, $L_{G^{(2)}}$, is uniformly bounded, because that of $\Gamma$, $L_\Gamma$, is,

$$|\mathcal{A}^{(2)}_{n,p}(x)| \leq \frac{L_\Gamma}{2\pi(\sigma n^{1/3})^{k-1}\delta^k} e^{-n\sigma^3\delta^3/6} e^{n^{1/3}\sigma(s_0 R_1 - (x\delta/2))}.$$

We deduce from this that for all $x$ in $[-s_0, \infty)$, $s_0 > 0$,

$$|\mathcal{A}^{(2)}_{n,p}(x)| \leq C(-s_0) e^{-x\varepsilon/2} e^{-n\sigma^3\delta^3/12},$$

when $n$ is large enough. Note also that $C(s_0)$ can be chosen to be a continuous nonincreasing function of $s_0$.

We now turn to $H^{(2)}_{\infty,k}$. Note that if $a \in \Gamma^{(2)}_\infty$, $a = te^{\pm i\pi/3}$, and $t \geq \delta\sigma n^{1/3}$. Recall that

$$H^{(2)}_{\infty,k}(x) = -\frac{1}{2\pi i}\int_{\Gamma^{(2)}_\infty} \frac{e^{-xa+a^3/3}}{a^k}\, da.$$

Hence,

$$|H^{(2)}_{\infty,k}(x)| \leq \frac{1}{2\pi}\int_{\Gamma^{(2)}_\infty} \frac{e^{-x\Re(a)+\Re(a^3)/3}}{|a|^k}|da| \leq \frac{1}{\pi}\int_{\delta\sigma n^{1/3}}^\infty \frac{e^{-xt/2-t^3/3}}{t^k}\, dt.$$

Now note that for $s_0 > 0$, $x \geq -s_0$ and $t > \varepsilon$, $e^{-xt/2} \leq e^{s_0 t/2 - x\varepsilon/2}$. We conclude from the last display that, when $n$ is large enough,

$$|\exp(-\varepsilon x) H^{(2)}_{\infty,k}(x)| \leq C(-s_0) e^{-\sigma^3\delta^3 n/6} e^{-x\varepsilon/2}$$

where again $C(y)$ can be chosen to be a nonincreasing continuous function of $y$.

As an aside, let us go back to the point we raised in the main text about $f$ not being defined when we cross the real axis. What we just did is to take the modulus of the quantity that appears inside the integral taken over $G^{(2)}$. When we worked on $\Re(f)$, we essentially focused on these quantities, since $\Re(\log(z)) = \log(|z|)$, when the log is defined. So the analysis we did for $\Re(f)$ applies to the situation when we first take the modulus of the quantity of interest, and hence we are rid of the problem created by the fact that the log is not defined when we cross the real axis at $R_1$.

B.2.4. *Behavior of the difference of our functions on $G^{(1)}$.* We first note that after changing variables through $a = \sigma n^{1/3}(z-c)$, $\Gamma_\infty^{(1)}$ is transformed into $G^{(1)}$. In other words, after doing this change of variables,

$$H_{\infty,k}^{(1)}(x) = -\frac{\sigma n^{1/3}}{2\pi i}\int_{G^{(1)}} e^{-\sigma n^{1/3}x(z-c)+n\sigma^3(z-c)^3/3}\frac{dz}{(n^{1/3}\sigma(z-c))^k}.$$

We also have $\exp(-\varepsilon x) = \exp(n^{1/3}\sigma(q-c)x)$, and therefore,

$$\exp(-\varepsilon x)H_{\infty,k}^{(1)}(x) = -\frac{\sigma n^{1/3}}{2\pi i}\int_{G^{(1)}} e^{-\sigma n^{1/3}x(z-q)+n\sigma^3(z-c)^3/3}\frac{dz}{(n^{1/3}\sigma(z-c))^k}.$$

Hence we have

$$\begin{aligned}&|\mathcal{A}_{n,p}^{(1)}(x) - \exp(-\varepsilon x)H_{\infty,k}^{(1)}(x)|\\ &\quad\le \frac{\sigma n^{1/3}}{2\pi}\int_{G^{(1)}} e^{-\sigma n^{1/3}x\Re(z-q)}\frac{|e^{n(f(z)-f(c))} - e^{n\sigma^3(z-c)^3/3}|}{|n^{1/3}\sigma(z-c)|^k}|dz|.\end{aligned}$$

- *The case $z \in \Gamma_0$.*

  Recall that if $z \in \Gamma_0$, $z - c = e^{i\theta}\varepsilon/(2\sigma n^{1/3})$, $\theta \in [\pi/3, 5\pi/3]$. We call

$$\mathfrak{I}_{\Gamma_0}(x) = \frac{\sigma n^{1/3}}{2\pi}\int_{G^{(1)}\cap\Gamma_0} e^{-\sigma n^{1/3}x\Re(z-q)}\frac{|e^{n(f(z)-f(c))} - e^{n\sigma^3(z-c)^3/3}|}{|n^{1/3}\sigma(z-c)|^k}|dz|.$$

Note that for $u, v \in \mathbb{C}$, $|e^u - e^v| \le \max(|e^u|, |e^v|)|u - v|$. This is easily seen if we write $\gamma(t) = v + (u-v)t$ and note that $e^u - e^v = \int_0^1 e^{\gamma(t)}\gamma'(t)\,dt$. Then, $|e^{\gamma(t)}| = \exp(\Re(v) + t\Re(u-v))$ and the result follows.

Hence, using the computations made in Section B.2.2,

$$\begin{aligned}&|e^{n(f(z)-f(c))} - e^{n\sigma^3(z-c)^3/3}|\\ &\quad\le \max(|e^{n(f(z)-f(c))}|, |e^{n\sigma^3(z-c)^3/3}|)n\left|f(z) - f(c) - \frac{\sigma^3}{3}(z-c)^3\right|\\ &\quad\le e^{n\sigma^3(2\Re((z-c)^3)+|z-c|^3)/6}n\Delta|z-c|^4.\end{aligned}$$

We also have $\sigma n^{1/3}(z-c) = e^{i\theta}\varepsilon/2$, and therefore

$$e^{n\sigma^3(2\Re((z-c)^3)+|z-c|^3)/6}n\Delta|z-c|^4 \le e^{\varepsilon^3/16}\frac{\Delta\varepsilon^4}{16\sigma^4 n^{1/3}}.$$

In other respects, $\Re(z-q)\sigma n^{1/3} = \sigma n^{1/3}(\Re(z-c) + \Re(c-q)) = \varepsilon(1+\cos(\theta)/2)$. We also note that the length of $\Gamma_0$ is $4\pi\varepsilon/(6\sigma n^{1/3})$. Therefore, we conclude that

$$\begin{aligned}&\int_{G^{(1)}\cap\Gamma_0} e^{-\sigma n^{1/3}x\Re(z-q)}\frac{|e^{n(f(z)-f(c))} - e^{n\sigma^3(z-c)^3/3}|}{|n^{1/3}\sigma(z-c)|^k}|dz|\\ &\quad\le C(-s_0)\exp(-\varepsilon x/2)\frac{4\pi\varepsilon}{6\sigma n^{1/3}}\left(\frac{\varepsilon}{2}\right)^{(-k)}e^{\varepsilon^3/16}\frac{\Delta\varepsilon^4}{16\sigma^4 n^{1/3}},\end{aligned}$$

where $C$ can be chosen to be a continuous nonincreasing function. In other words, for $x \in [-s_0, \infty)$, when $n$ is large enough,

$$\mathfrak{I}_{\Gamma_0}(x) \leq \frac{C(-s_0)\exp(-x\varepsilon/2)}{n^{1/3}}.$$

- *The case $z \in \Gamma_1$.*

  When $z \in \Gamma_1 \cap G^{(1)}$, $z = c + te^{\pm i\pi/3}$, $\varepsilon/(2\sigma n^{1/3}) \leq t \leq \delta$.
  We call

$$\mathfrak{I}_{\Gamma_1}(x) = \frac{\sigma n^{1/3}}{2\pi}\int_{G^{(1)}\cap\Gamma_1} e^{-\sigma n^{1/3}x\Re(z-q)}\frac{|e^{n(f(z)-f(c))} - e^{n\sigma^3(z-c)^3/3}|}{|n^{1/3}\sigma(z-c)|^k}|dz|.$$

  Going through the same steps as before, we find that

$$\begin{aligned}|e^{n(f(z)-f(c))} - e^{n\sigma^3(z-c)^3/3}| &\leq e^{n\sigma^3(2\Re((z-c)^3)+|z-c|^3)/6} n\Delta|z-c|^4\\ &\leq e^{-n\sigma^3t^3/6}\Delta nt^4.\end{aligned}$$

In other respects, $\sigma n^{1/3}\Re(z-q) = tn^{1/3}\sigma/2 + \varepsilon$. Therefore, for $x \in [-s_0, \infty)$, with $s_0 > 0$, we have $-x\sigma n^{1/3}\Re(z-q) \leq s_0(tn^{1/3}\sigma/2 + \varepsilon) - \varepsilon x$. Hence,

$$\mathfrak{I}_{\Gamma_1}(x) \leq \frac{\sigma n^{1/3}}{2\pi}e^{s_0\varepsilon}e^{-\varepsilon x}\int_{\varepsilon/(2\sigma n^{1/3})}^{\delta} e^{s_0 tn^{1/3}\sigma/2}\frac{e^{-n\sigma^3t^3/6}\Delta nt^4}{(n^{1/3}\sigma t)^k}\,dt.$$

After changing variables to $v = \sigma n^{1/3}t$, we get

$$\mathfrak{I}_{\Gamma_1}(x) \leq \frac{\Delta e^{s_0\varepsilon}e^{-\varepsilon x}}{2\pi\sigma^4 n^{1/3}}\int_{\varepsilon/2}^{\infty} e^{s_0 v/2}e^{-v^3/6}v^{4-k}\,dv.$$

Hence, here again, we can find a continuous nondecreasing function $C$ such that for $x \geq -s_0$, $s_0 > 0$ and for $n$ large enough,

$$\mathfrak{I}_{\Gamma_1}(x) \leq \frac{C(-s_0)e^{-x\varepsilon/2}}{n^{1/3}}.$$

B.2.5. *Conclusion.* The expression $|\kappa_{n,p}A_{n,p} - e^{-\varepsilon x}H_{\infty,k}(x)|$ was our initial center of interest. We have the simple bound

$$\begin{aligned}&|\kappa_{n,p}A_{n,p} - e^{-\varepsilon x}H_{\infty,k}(x)|\\ &\quad\leq |\mathcal{A}^{(1)}_{n,p}(x) - \exp(-\varepsilon x)H^{(1)}_{\infty,k}(x)| + |\mathcal{A}^{(2)}_{n,p}(x)| + \exp(-\varepsilon x)|H^{(2)}_{\infty,k}(x)|\\ &\quad\leq \mathfrak{I}_{\Gamma_0}(x) + \mathfrak{I}_{\Gamma_1}(x) + |\mathcal{A}^{(2)}_{n,p}(x)| + |\exp(-\varepsilon x)H^{(2)}_{\infty,k}(x)|\\ &\quad\leq \frac{C(-s_0)e^{-\varepsilon x/2}}{n^{1/3}},\end{aligned}$$

for $C$ a nonincreasing continuous function. This bound is valid if $x \in [-s_0, \infty)$, $s_0 > 0$ and when $n$ is large enough.

So Lemma B.1 is shown.

**B.3. Convergence of $B_{n,p}$.** We again work in the general case where there is a root of multiplicity $k$ at $c(\Sigma_p, n, p)$. With notation similar to the ones above, we have

$$B_{\widetilde{\Sigma}_{p+k},n,p+k}(s) = B_{n,p}(s) = \frac{n^{1/3}\sigma}{2\pi i}\int_{\Xi} e^{n^{1/3}\sigma(z-q)}e^{-nf(z)}\frac{(c-z)^k}{c^k}\,dz.$$

Hence,

$$\mathcal{B}_{n,p}(x) \triangleq B_{n,p}(s)/\kappa_{n,p} = \frac{(\sigma n^{1/3})^{k+1}}{2\pi i}\int_{\Xi} e^{n^{1/3}\sigma s(z-q)}e^{-n(f(z)-f(c))}(z-c)^k\,dz.$$

The aim of this subsection is to show the following lemma.

Lemma B.2. *Let $f$ satisfy Lemma 2. Suppose*

$$\exists \delta > 0 : \forall s \qquad |s-c| < \delta \Rightarrow \left|\frac{f^{(4)}(s)}{4!}\right|\delta \leq \frac{\sigma^3}{6}.$$

*Then*

$$\begin{aligned}&\forall s_0 \in \mathbb{R}, \exists C(s_0), \exists N_0(s_0) : s > s_0, n > N_0\\ &\qquad \Rightarrow |B_{n,p}(s)/\kappa_{n,p} - e^{\varepsilon s}J_{\infty,k}(s)| \leq \frac{C(s_0)\exp(-\varepsilon s/2)}{n^{1/3}}.\end{aligned}$$

We now turn to proving Lemma B.2.

B.3.1. *Notation and preliminary computations.* Recall that we denote by $\mathcal{C}$ the circle of center $c$ and radius $\delta$. We call $\mathcal{D}$ the corresponding disc. We split $\Xi$ into $\Xi = X^{(1)} \cup X^{(2)}$, where $X^{(1)}$ is the part of $\Xi$ that is inside $\mathcal{D}$. Note that under our assumptions about $\delta$ and $e$, the intersection of $\Xi$ and $\mathcal{C}$ is on $\Xi_1 \cup \overline{\Xi_1}$, that is, on a section of $\Xi$ where this contour is parametrized as $c + te^{\pm 2i\pi/3}$, $t \in \mathbb{R}$.

We call $\Xi_\infty^{(1)}$ the image of $X^{(1)}$ under the map $z \mapsto \sigma n^{1/3}(z-c)$. Of course, everything has been done so that this is a subset of $\Xi_\infty$. Let us denote $\Xi_\infty^{(2)} = \Xi_\infty \setminus \Xi_\infty^{(1)}$.

Let us call $\mathcal{B}_{n,p}(x) = B_{n,p}(x)/\kappa_{n,p}$ and $\mathcal{B}_{n,p}(x) = \mathcal{B}_{n,p}^{(1)}(x) + \mathcal{B}_{n,p}^{(2)}(x)$, where the subscript indicates that $\mathcal{B}_{n,p}^{(i)}(x)$ is the contribution of the integral defining $\mathcal{B}_{n,p}(x)$ over $X^{(i)}$.

We similarly split $J_{\infty,k}$ into $J_{\infty,k} = J_{\infty,k}^{(1)} + J_{\infty,k}^{(2)}$ where now the subscripts refer to the contribution of the integrals over $\Xi_\infty^{(1)}$ and $\Xi_\infty^{(2)}$.

Note that using the same arguments as before, we have, inside $\mathcal{D}$,

$$\Re(f(c) - f(z)) \leq \frac{\sigma^3}{6}(-2\Re((z-c)^3) + |z-c|^3).$$

So for $z = c + te^{\pm i2\pi/3}$, $\Re(f(c) - f(z)) \leq -\sigma^3 t^3/6$. Also, by arguments similar to the ones we used before, we have

$$\max_{z\in\Xi_\infty^{(2)}} \Re(-f(z)) \leq \Re(f(c+\delta e^{2i\pi/3})) \leq -f(c) - \frac{\sigma^3}{6}\delta^3.$$

B.3.2. *Behavior of our functions on* $X^{(2)}$ *and* $\Xi_\infty^{(2)}$. We first focus on $\mathcal{B}_{n,p}^{(2)}(x)$. Note that for $z \in X^{(2)}$, we have $|z - c| \leq \sqrt{3c^2 + (R_2 + c)^2}$ and $-(R_2 + c) \leq \Re(z - q) \leq -\varepsilon/(2\sigma n^{1/3})$. The second inequality comes from the fact that $\Xi_0$ is a circle of radius $3\varepsilon/(\sigma n^{1/3})$. Now we have

$$|\mathcal{B}_{n,p}^{(2)}(x)| \leq \frac{(\sigma n^{1/3})^{k+1}}{2\pi}\int_{\Xi_\infty^{(2)}} e^{n^{1/3}\sigma x\Re(z-q)} e^{n\Re(f(c)-f(z))}|z-c|^k|dz|.$$

Because for $x$ in $[-s_0, \infty)$, $x\Re(z - q) \leq -x\varepsilon/(2\sigma n^{1/3}) + s_0(R_2 + c)$, we conclude using the results put forth in the previous subsection, that

$$|\mathcal{B}_{n,p}^{(2)}(x)| \leq C(-s_0)e^{-\sigma^3\delta^3 n/12}e^{-\varepsilon x/2}.$$

On the other hand, since $\Xi_\infty^{(2)}$ is parametrized as $z = te^{\pm i2\pi/3}$, with $\delta\sigma n^{1/3} \leq t < \infty$, we obtain along the lines of the proof done for $H_{\infty,k}^{(2)}$ that

$$J_{\infty,k}^{(2)}(x) \leq e^{-3\varepsilon x/2}e^{-\sigma^3\delta^3 n/6}C(s_0).$$

B.3.3. *Behavior of the difference of our functions on* $X^{(1)}$. Here again, by the change of variables $a = n^{1/3}\sigma(z - c)$, $\Xi_\infty^{(1)}$ is mapped to $X^{(1)}$. Hence we can write

$$\begin{aligned}
&|\mathcal{B}_{n,p}^{(1)}(x) - e^{\varepsilon x}J_{\infty,k}^{(1)}(x)| \\
&\quad \leq \frac{\sigma n^{1/3}}{2\pi}\int_{X^{(1)}} e^{\sigma n^{1/3}x\Re(z-q)}|\sigma n^{1/3}(z-q)|^k \\
&\qquad\qquad \times |e^{n(f(c)-f(z))} - e^{-n\sigma^3(z-c)^3/3}||dz|.
\end{aligned}$$

Splitting the problem into first $\Xi_0$ and then $\Xi_1$, and repeating the approach used in the study of $\mathcal{A}^{(1)}$ together with the new estimates of $f(c) - f(z)$, we get

$$|\mathcal{B}_{n,p}^{(1)}(x) - e^{\varepsilon x}J_{\infty,k}^{(1)}(x)| \leq C(-s_0)e^{-\varepsilon x/2}/n^{1/3}.$$

The combination of this bound and those for $J_{\infty,k}^{(2)}$ and $\mathcal{B}_{n,p}^{(2)}(x)$ shows Lemma B.2.

**Acknowledgments.** I am very grateful to Professor Iain Johnstone for numerous discussions and to Professor Persi Diaconis for discussions and useful references. I would also like to thank Professor Donald St. P. Richards for correspondence and references, Professor Craig A. Tracy for references and Professor Alice Whittemore for her support. I am grateful to an anonymous referee for constructive comments and useful references.

## REFERENCES

[1] Anderson, T. W. (1963). Asymptotic theory for principal component analysis. *Ann. Math. Statist.* **34** 122–148. MR0145620

[2] Anderson, T. W. (2003). *An Introduction to Multivariate Statistical Analysis*, 3rd ed. Wiley, Hoboken, NJ. MR1990662

[3] Bai, Z. D. (1999). Methodologies in spectral analysis of large-dimensional random matrices, a review. *Statist. Sinica* **9** 611–677. MR1711663

[4] Bai, Z. D. and Silverstein, J. W. (1998). No eigenvalues outside the support of the limiting spectral distribution of large-dimensional sample covariance matrices. *Ann. Probab.* **26** 316–345. MR1617051

[5] Baik, J. (2006). Painleve formulas of the limiting distributions for nonnull complex sample covariance matrices. *Duke Math. J.* **133** 205–235. MR2225691

[6] Baik, J., Ben Arous, G. and Péché, S. (2005). Phase transition of the largest eigenvalue for non-null complex sample covariance matrices. *Ann. Probab.* **33** 1643–1697. MR2165575

[7] Borodin, A. (1999). Biorthogonal ensembles. *Nuclear Phys. B* **536** 704–732. MR1663328

[8] Böttcher, A. and Silbermann, B. (1999). *Introduction to Large Truncated Toeplitz Matrices*. Springer, New York. MR1724795

[9] Desrosiers, P. and Forrester, P. J. (2006). Asymptotic correlations for Gaussian and Wishart matrices with external source. *Int. Math. Res. Not.* **2006** Article ID 27395. MR2250019

[10] Dieng, M. (2005). Distribution functions for edge eigenvalues in orthogonal and symplectic ensembles: Painlevé representations. *Int. Math. Res. Not.* **2005** 2263–2287. MR2181265

[11] El Karoui, N. (2006). A rate of convergence result for the largest eigenvalue of complex white Wishart matrices. *Ann. Probab.* To appear.

[12] El Karoui, N. (2003). On the largest eigenvalue of Wishart matrices with identity covariance when $n, p$ and $p/n \to \infty$. Available at arXiv:math.ST/0309355.

[13] El Karoui, N. (2004). New results about random covariance matrices and statistical applications. Ph.D. dissertation, Stanford Univ.

[14] Forrester, P. J. (1993). The spectrum edge of random matrix ensembles. *Nuclear Phys. B* **402** 709–728. MR1236195

[15] Forrester, P. J. (2006). Eigenvalue distributions for some correlated complex sample covariance matrices. Available at arxiv:math-ph/0602001.

[16] Gohberg, I., Goldberg, S. and Krupnik, N. (2000). *Traces and Determinants of Linear Operators*. Birkhäuser, Basel. MR1744872

[17] Gravner, J., Tracy, C. A. and Widom, H. (2001). Limit theorems for height fluctuations in a class of discrete space and time growth models. *J. Statist. Phys.* **102** 1085–1132. MR1830441

[18] Gray, R. M. (2006). Toeplitz and circulant matrices: A review. *Foundations and Trends in Communications and Information Theory* **2** 155–239. Available at http://ee.stanford.edu/~gray/toeplitz.pdf.

[19] Grenander, U. and Szegö, G. (1958). *Toeplitz Forms and Their Applications*. Univ. California Press, Berkeley. MR0094840

[20] Gross, K. and Richards, D. (1989). Total positivity, spherical series, and hypergeometric functions of matrix argument. *J. Approx. Theory* **59** 224–246. MR1022118

[21] Guionnet, A. and Zeitouni, O. (2000). Concentration of the spectral measure for large matrices. *Electron. Comm. Probab.* **5** 119–136. MR1781846

[22] Horn, R. and Johnson, C. (1990). *Matrix Analysis*. Cambridge Univ. Press. MR1084815

[23] James, A. T. (1964). Distributions of matrix variates and latent roots derived from normal samples. *Ann. Math. Statist.* **35** 475–501. MR0181057

[24] Johansson, K. (2000). Shape fluctuations and random matrices. *Comm. Math. Phys.* **209** 437–476. MR1737991

[25] Johnstone, I. (2001). On the distribution of the largest eigenvalue in principal component analysis. *Ann. Statist.* **29** 295–327. MR1863961

[26] Ledoux, M. (2001). *The Concentration of Measure Phenomenon*. Amer. Math. Soc., Providence, RI. MR1849347

[27] Marčenko, V. A. and Pastur, L. A. (1967). Distribution of eigenvalues in certain sets of random matrices. *Mat. Sb. (N.S.)* **72** 507–536. MR0208649

[28] Olver, F. W. J. (1974). *Asymptotics and Special Functions*. Academic Press, New York–London. MR0435697

[29] Paul, D. (2007). Asymptotics of sample eigenstructure for a large dimensional spiked covariance model. *Statist. Sinica*. To appear.

[30] Reed, M. and Simon, B. (1972). *Methods of Modern Mathematical Physics. I. Functional Analysis*. Academic Press, New York. MR0751959

[31] Silverstein, J. W. and Choi, S.-I. (1995). Analysis of the limiting spectral distribution of large-dimensional random matrices. *J. Multivariate Anal.* **54** 295–309. MR1345541

[32] Simon, S., Moustakas, A. and Marinelli, L. (2005). Capacity and character expansions: Moment generating function and other exact results for mimo correlated channels. Available at arxiv:cs.IT/0509080.

[33] Soshnikov, A. (2000). Determinantal random point fields. *Russian Math. Surveys* **55** 923–975. MR1799012

[34] Soshnikov, A. (2002). A note on universality of the distribution of the largest eigenvalues in certain sample covariance matrices. *J. Statist. Phys.* **108** 1033–1056. MR1933444

[35] Tracy, C. and Widom, H. (1994). Level-spacing distribution and the Airy kernel. *Comm. Math. Phys.* **159** 151–174. MR1257246

[36] Tracy, C. and Widom, H. (1996). On orthogonal and symplectic matrix ensembles. *Comm. Math. Phys.* **177** 727–754. MR1385083

[37] Tracy, C. and Widom, H. (1998). Correlation functions, cluster functions and spacing distributions for random matrices. *J. Statist. Phys.* **92** 809–835. MR1657844

[38] Tulino, A. and Verdú, S. (2004). *Random Matrix Theory and Wireless Communications*. *Foundations and Trends in Communications and Information Theory* **1**. Now Publishers, Hanover, MA.

[39] van der Vaart, A. W. (1998). *Asymptotic Statistics*. Cambridge Univ. Press. MR1652247

[40] Wachter, K. W. (1978). The strong limits of random matrix spectra for sample matrices of independent elements. *Ann. Probab.* **6** 1–18. MR0467894

[41] Widom, H. (1999). On the relation between orthogonal, symplectic and unitary matrix ensembles. *J. Statist. Phys.* **94** 347–363. MR1675356

Department of Statistics
University of California
367 Evans Hall
Berkeley, California 94720
USA
E-mail: nkaroui@stat.berkeley.edu